\documentclass{amsart}

\usepackage{amsfonts,amsthm,amsmath,amsfonts,latexsym,amssymb}
\usepackage[all]{xy}
\usepackage{graphicx,upref,calligra}
\usepackage[svgnames]{xcolor}
\usepackage{mathrsfs}
\usepackage[T1]{fontenc}
\usepackage{float}
\usepackage{bbm}
\numberwithin{equation}{section}

\newtheorem{fed}{Definition}[section]
\newtheorem{teo}[fed]{Theorem}
\newtheorem{lem}[fed]{Lemma}
\newtheorem{cor}[fed]{Corollary}
\newtheorem{pro}[fed]{Proposition}
\newtheorem{rem}[fed]{Remark}

\newtheorem{exa}[fed]{Example}

\theoremstyle{definition}

\theoremstyle{remark}

\definecolor{azul}{rgb}{0.1,0.6,0.86}
\definecolor{titleblue}{rgb}{0.13,0.49,0.69}
\definecolor{mylred}{rgb}{0.85,0.24,0.2}
\definecolor{myblue}{rgb}{0,0.33,0.55}
\definecolor{myyellow}{rgb}{0.42,0.24,0.52}
\definecolor{mygreen}{rgb}{0.12,0.5,0.29}
\definecolor{myred}{rgb}{0.74,0.13,0.13}
\definecolor{mylblue}{rgb}{0.2,0.75,1}
\definecolor{mylgreen}{rgb}{0.68,0.98,0.6}
\definecolor{mylyellow}{rgb}{0.86,0.85,0.55}
\definecolor{myllyellow}{rgb}{0.87,0.86,0.56}
\definecolor{naranja}{RGB}{249,153,96}
\definecolor{sidebardarkcolor}{rgb}{0.21,0.31,0.40}
\definecolor{sidebarlightcolor}{rgb}{0.7,0.77,0.836}

\def\noi{\noindent}

\def\EOE{\hfill $\blacktriangle$}

\def\bdem{\begin{proof}}
\def\edem{\end{proof}}

\def\bm{\left(\begin{array}}
\def\em{\end{array}\right)}
\def\ben{\begin{enumerate}}
\def\een{\end{enumerate}}
\def\barr{\begin{array}}
\def\earr{\end{array}}

\def\eps{\varepsilon}
\def\fii{\varphi }
\def\la{\lambda}
\def\La{\Lambda}

\def\w{\omega}
\def\W{\Omega}
\def\si{\sigma}

\def\N{\mathbb{N}}
\def\Z{\mathbb{Z}}
\def\R{\mathbb{R}}
\def\C{\mathbb{C}}

\def\D{\mathbb{D}}
\def\T{\mathbb{T}}

\def\Q{\mathbb{Q}}
\def\K{\mathbb{K}}

\def\cD{\mathcal{D}}
\def\ede{\mathcal{D}}

\def\G{\mathcal{G}}

\def\ese{\mathcal{S}}

\def\eme{\mathcal{M}}

\def\ewe{\mathcal{W}}

\def\a{\mathbf a}

\newcommand{\peso}[1]{ \quad \text{ \rm  #1 } \quad }

\newcommand{\sub}[2]{{#1}_{\mbox{\tiny{${#2}$}}}}

\DeclareMathOperator{\leqm}{\preccurlyeq}

\newcommand{\pint}[1]{\displaystyle \left \langle\, #1 \, \right\rangle}

\newcommand{\swb}{\mathcal{S}(G)}
\newcommand{\swbd}{\mathcal{S}(\widehat{G})}
\newcommand{\swbp}[1]{\mathcal{S}(#1)}
\newcommand{\ag}{\mathcal{A}(G)}
\newcommand{\agd}{\mathcal{A}(\widehat{G})}

\newcommand{\m}[1]{|#1|}

\def\R{\mathbb R} 
\def\Z{\mathbb Z}
\def\C{\mathbb C}
\def\Q{\mathbb Q}
\def\T{\mathbb T}
\def\N{\mathbb N} 

\def\L{\Lambda}
\def\Chi{\raise .3ex \hbox{\large $\chi$}}

\def\a{\alpha}

\def\w{\omega}

\def\cD{\mathcal{D}}

\def\G{{\widehat{G}}}

\begin{document}
% ----------------------------------------------------------------
% ----------------------------------------------------------------

\title[Quasicrystals in LCA groups]{Existence of quasicrystals and universal stable sampling and interpolation in LCA groups}

\author{Elona Agora}
\address{Instituto Argentino de Matemática "Alberto P. Calderón" (IAM-CONICET), Buenos Aires, Argentina}
\email{elona.agora@gmail.com}
\thanks{The first author was supported in part by Grants:  MTM2016-75196-P (MINECO / FEDER, UE),  PIP 112201501003553CO,  UBACyT 20020130100422BA, PICT 2014-1480 (ANPCyT)}

\author{Jorge Antezana}
\address{Departamento de Matemática, Universidad Nacional de La Plata and, 
 Instituto Argentino de Matemática "Alberto P. Calderón" (IAM-CONICET), Buenos Aires, Argentina}
\email{antezana@mate.unlp.edu.ar}
\thanks{The second author was supported in part  by Grants: CONICET-PIP 152,  UNLP-11X585, MTM2016-75196-P}

\author{Carlos Cabrelli }
\address{Departamento de Matemática, Universidad de Buenos Aires and, 
Instituto de Matemática "Luis Santaló" (IMAS-CONICET-UBA), Buenos Aires, Argentina}
\email{cabrelli@dm.uba.ar}
\thanks{The third author was supported in part by Grants:  PICT 2014-1480 (ANPCyT), CONICET PIP 11220110101018,
UBACyT 20020130100403BA, UBACyT 20020130100422BA}

\author{Basarab Matei}
\address{Institut Galilée and, Université Paris 13, Paris, France}
\email{matei@lipn.univ-paris13.fr}

\subjclass[2010]{Primary 42C15, 94A20; Secondary 42C30, 43A25}

\keywords{ Quasicrystals, Universal Sampling and Interpolation, Landau-Beurling's densities,  Poisson measures,  Locally compact abelian groups}

\begin{abstract}
We characterize all the locally compact abelian (LCA)  groups that contain quasicrystals (a class of model sets). Moreover, we describe all possible quasicrystals in the group constructing an appropriate lattice associated with the cut and project scheme that produces it. On the other hand, if an LCA group $G$ admits a simple quasicrystal, we prove that recent results of Meyer and Matei for the case of the Euclidean space $\R^n$ can be extended to $G$. More precisely, we prove that  simple quasicrystals are universal sets of stable sampling and universal sets of stable interpolation in generalized Paley-Wiener spaces. 
\end{abstract}

\maketitle

\section{Introduction} \label{Intro}

Quasicrystals  are non-periodic structures discovered by Shechtman  in 1984, while studying materials 
whose $X$-ray diffractions spectra present such non-periodic behaviors (see \cite{Lag},~\cite {Sh}). 
One of the best mathematical models for quasicrystals are model sets introduced  by Meyer in \cite{Meyer2} many years before
 (see \cite{AM}, \cite{Lag1} and \cite{Moody2} for a survey on model sets and quasicrystals). 

\medskip

In fact, Meyer proposed  a method called cut and project to construct model sets. Given a locally compact  group $G$ and $m \in \N$, a \textbf{cut and project scheme}  (CP-scheme)  for $G$ is a triple $(\R^m,G,H)$ where $H\subseteq \R^m \times G$ is a lattice such that 
  the canonical projections, $p_1$ and $p_2$  satisfy: $ p_1:H \rightarrow \R^m$ is one to one and 
 $p_2: H \rightarrow G$  has  dense range.   
 We will call a CP-scheme \textbf{complete} if, in addition,  $p_1: H \rightarrow \R^m$  has  dense range and $p_2:H \rightarrow G$ is one to one. 

\medskip

Given a complete CP-scheme  $(\R^m,G,H)$ and a Riemann integrable set (i.e. the Lebesgue measure of its boundary is zero) with non-empty interior $S\subset\R^m$, the \textbf{model set} associated to this scheme and the set $S$, is defined by
\begin{equation}
\label{(3.1)}
\L_S=\{p_2(h): h \in H, \,p_1(h)\in S\}.
\end{equation} 

Throughout this paper,  following \cite{MM1}, these model sets will be called \textbf{quasicrystals}.
When $m= 1$ and $S \subset \R$ is an interval the quasicrystal will be called {\bf simple}.

\subsection*{Quasicrystals and Fourier analysis} 

In \cite{MM1}  Matei and Meyer  discovered that quasicrystals play an important role in Fourier Analysis, more precisely in the theory of sampling and interpolation of band limited functions (see also \cite{MM2}). Before describing their results, let us recall the basic definitions. 

\medskip

Let $G$ be a locally compact abelian (LCA) group, and let $K$ be a compact subset of $\G$, the dual group of $G$. Recall that the \textbf{Paley-Wiener space} $PW_K$ consists on all square integrable functions defined on $G$ whose Fourier transform vanishes (almost everywhere) outside $K$. For this space, a set $\Lambda\subseteq G$ is a \textbf{stable sampling set} if there exist constants $A, B>0$ such that for any $f\in PW_K$,
$$
A\|f\|_2^2 \leq \sum _{\la \in \La} |f(\la)|^2 \leq B\|f\|^2_2.
$$

\medskip

If $\La$ is  a stable sampling set, then for any $f\in PW_K$, the values $\{f(\la)\}_{\la\in\La}$ contain enough information to recover completely $f$. On the other hand, a set  $\Gamma$ is an \textbf{ stable interpolation set} for $PW_K$ if  for every  $\{c_\gamma \}_{\gamma \in \Gamma}  \in \ell^2(\Gamma)$,
 the interpolation problem
$$
f(\gamma)=c_\gamma, 
$$
has a solution $f\in PW_K$. 

\medskip

Intuitively, a sampling  set needs to have enough information (in terms of number of points) to be able to recover the sampled function. This means that the more points we have, the more information we obtain. On the other hand, each point in an interpolation set imposes an extra restriction to find an interpolation function. Thus, the less points we have, the better it is. This intuition is formalized by means of the so called \textbf{Beurling-Landau densities}. If $\La\subset \R^d$ the \textbf{upper and lower Landau-Beurling densities} are  defined by
$$
\cD^+ (\La)= \lim_{\ell \to \infty} \sup_{x\in \R^d} \frac{  \# ( \La \cap Q_\ell (x) ) }{\ell^d} \quad \text{and} \quad \cD^- (\La)= \lim_{\ell \to \infty} 
\inf_{x\in \R^d} \frac{ \# (\La \cap  Q_\ell (x))} {\ell^d}
 $$ 
respectively. Here,  $ Q_\ell (x) $ denotes the cube centered at $x$, of side length $\ell$.  If $\cD^+=\cD^-$, we say that the set has \textbf{uniform density} and in this case is simply denoted by $\cD$.  Landau proved the following necessary conditions in \cite{Landau}: 
\begin{itemize}
\item [(i)] A sampling set $\La$ for $PW_{\Omega}$ satisfies  $\cD^-(\La)\geq |\Omega|$;
\item [(ii)] An interpolation set $\La$ for $PW_{\Omega}$ satisfies  $\cD^+(\La)\leq |\Omega|$. 
\end{itemize}
Later on, this result was extended to LCA groups by Gr\"ochenig, Kutyniok and Seip \cite{neclca}.  In this case, the Lebesgue measure is replaced by the Haar measure of $\widehat{G}$ and, $\cD^+$ and  $\cD^-$ denote the upper and lower Landau-Beurling densities generalized to the group setting \cite{neclca} (see Section \ref{LBdensities} for precise definitions).
Roughly speaking, the Landau-Beurling's densities compare the distribution of the points of $\La$ with that of a reference set,  
which for instance in the case of $\R^d$ is the lattice $\Z^d$. 

\medskip

Beurling showed in \cite{Beurling} the following kind of converse of Landau's result, if $\Omega$ is an interval of the real line.

\begin{teo}\label{TeoBeurling}  
Let $I\subseteq \R$ be a bounded interval and $\La\subseteq \R$  a uniformly discrete sequence (i.e. there exists $M>0$ so that $|\la_1-\la_2|>M$ for all $\la_1,\la_2\in \La$ different). Then  
\begin{itemize}
\item[(i)]  If $\cD^-(\La)>|I|$ then $\La$  is a stable sampling set for $\text{PW}_I$;
\item[(ii)] If $\cD^+(\La)<|I|$ then $\La$ is an stable interpolation set for $\text{PW}_I$. 
\end{itemize}
\end{teo}

Just by considering $\Omega=[0,1/2-\eps)\cup[1,3/2-\eps)$ for any small value of $\eps>0$, and $\La=\Z$, it is easy to see that a similar result is no longer true for more general sets. However, Olevskii and Ulanovskii found sets $\La\subset \R$  with uniform density and the remarkable property of being stable sampling sets (resp. interpolating sets) for any $PW_{\Omega}$,  such that $|\Omega|<D(\La)$ (resp. $|\Omega|>D(\La)$ and $\Omega$ Riemann integrable) (see \cite{Olev2} and \cite{Olev1}). This achievement is surprising because no assumption on the structure of $\Omega$ is required. Such sets have been called  \textbf{universal stable sampling set} (resp. \textbf{universal stable interpolation set}). 

\medskip

In \cite{MM1} Matei and Meyer proved the following extension of Beurling theorem.

\begin{teo}
Let $\L_I \subset \R^d$ be a simple quasicrystal  and let $K  \subset \R^d $ be a compact set. If $|K|$ denotes the standard Lebesgue measure of $K$ then
\begin{itemize}
\item[(i)] If $\cD(\L_I) >  |K|$ then $\L_I$ is a  stable sampling set for $PW_K$; 
\item[(ii)] If $K$ is Riemann integrable and $\cD(\L_I) < |K|$ then $\L_I$ is an stable interpolation set for $PW_K$. 
\end{itemize}
\end{teo}

In particular, this proves that simple quasicrystals in $\R^d$ are universal sets of sampling (resp. interpolation). This gave a very simple way to construct such  sets. Since then, quasicrystals played a key role in several recent advances in Fourier Analysis (see for instance \cite{GeL} and \cite{LevOlev}).

\subsection*{LCA group setting and main results}

Due to the importance of quasicrystals in Fourier Analysis,  a natural question is whether or not any LCA group has quasicrystals. In the case it has, the next natural question is whether or not they are universal sets of stable sampling and interpolation. 

\medskip

Let $G$ be a given LCA group such that its dual $\G$ is metrizable and compactly generated.  These are the standard assumptions to study  sampling and interpolation problems in LCA groups (see  \cite{neclca}). Then, by the structure theorems, $G$ is isomorphic to $\R^d \times \T^m \times \D$, where $\D$ is a countable discrete group (see for instance \cite{DE}).  To avoid some pathological cases, we will assume that $G$ is neither compact nor discrete.

\medskip

With respect to the problem of the existence of quasicrystals, if $G$ is isomorphic to the Lie group $\R^d \times \T^m$ then it is not difficult to see that the answer is positive. However, when the discrete part $\D$ appears, the problem becomes much more complicated.  This complication comes from the torsion of $\D$ (see Subsection \ref{initials} for a detailed discussion on this). Given $m\in\N$, our first main result is a complete characterization of those groups $G$ such that there exists a complete CP-scheme $(\R^m, G, H)$. In particular, this solves the problem of the existence of quasicrystals. More precisely we prove the following theorem, which is the main result of this work.

\begin{teo}\label{imposible} 
Let $m\in\N$ and let $G=\R^d\times \T^\ell\times \D$. Then, there exists a complete CP-scheme $(\R^m, G, H)$ if and only if $\D$ does not have a copy of $\Z_p^{m+d+1}$ for any prime $p$.
\end{teo} 
For the groups that do contain quasicrystals, we construct  essentially all the quasicrystals in the group describing all possible lattices involved in  the cut and 
project scheme  that produce the quasicrystal.

\medskip

Once we know which groups contains quasicrystals, we study the problem of universal sets of stable sampling and stable interpolation. Since the group $G$ is isomorphic to $\R^d\times \T^\ell\times \D$, the idea is to find a universal set of stable sampling (resp. interpolation) in each connected component 
$\R^d\times \T^\ell\times\{d\}$. However, the choice of these universal sets is a difficult problem. For instance, the simple idea of taking the same set in each component leads to sampling sets which are not universal (see the begining of Section \ref{DSI} for more details). So, in order to get universal sets of stable sampling or stable interpolation we should consider a much more subtle combination of such sets in each component. To solve this problem we use  quasicrystals. Our approach essentially follows the line of proof of \cite{MM1} adapted to our setting. Firstly, we extend the following formula to our setting.
\begin{teo}\label{nl} 
Let $H$  be a lattice on $ \R^m \times G$ such that $\left.p_1\right|_H$ and $\left.p_2\right|_H$ are one to one, and $p_1(H)$ is dense in $\R^m$. Let $\fii$ in the Schwartz class of $\R^m$ and $\psi$ in the Schwartz-Bruhat class of $G$. Therefore
$$
\lim_{r\to \infty}\,  \frac{1}{r^m} \,  \sum_{h \in H} \fii\left(\frac{p_1(h)-a}{r}\right)\psi(p_2(h))
\, = \, \frac{1}{\m{H}}\ \int_{\R^m}{\fii}(x)\,dx\,\int_G {\psi}(g)\,dm_G
$$
uniformly in $a\in \R^m$. 
\end{teo}
This result is interesting in its own because of its connection with Poisson measures and Poisson summation formul{\ae}  (see Subsection \ref{PM}).

\medskip

Let  $\Gamma \subseteq \R^m \times \G $ be the dual lattice of $H$ and consider the canonical projections
$$
q_1: \R^m\times \G \to \R^m \peso {and} q_2: \R^m\times \G \to \G.
$$
It can be seen that 
the projections $ \left.q_1\right|_\Gamma \,$ and $ \,\left.q_2\right|_\Gamma$ are injective and have dense range on  
$\R^m$ and $\G$ respectively (see Lemma 2, pg. 41, \cite{Meyer2}). Given a compact set  $K\subset \G$,  we define the set
\begin{equation} \label{dualm}
M_K := \{q_1(\gamma): \gamma\in \Gamma, \, q_2(\gamma) \in K \} \subseteq \R^m.
\end{equation}

\medskip

Now,  using Theorem \ref{nl}, the next  duality theorem follows using the same steps as in the case of $\R^d$.

\begin{teo}[Duality] \label{duality} 
Let $S$ be a compact set of $\R^m$ and let $K$ be a compact Riemann integrable
set of $\G$.  Then
\begin{itemize}
\item[(i)] If  $M_K$ is an stable interpolation set for $PW_S$, then $\La_S$ is a stable sampling set for $PW_K$;
\item[(ii)]  If $M_K$ is a stable sampling set for $PW_{\tilde{S}}$, then  $\La_{\tilde{S}}$ is a set of stable interpolation for 
$PW_K$, where $\tilde{S}$ is a slight dilation of $S$.
\end{itemize} 
\end{teo}

Finally, as a consequence of this duality theorem, the classical Beurling theorem \cite{Beurling}, 
and the extension to groups of Landau's theorem  \cite{neclca}, we obtain a Beurling-type theorem for quasicrystals  in LCA groups. 

\begin{teo}\label{thsuf} 
Let $\L_I \subset G$ be a simple quasicrystal  and let $K  \subset \G $ be a compact 
set. 
\begin{itemize}
\item[(i)] If $\cD(\L_I) >  \mu_{\G}(K)$ then $\L_I$ is a  stable sampling set for $PW_K$. 
\item[(ii)] If $K$ is Riemann integrable and $\cD(\L_I) < \mu_{\G}(K)$ then $\L_I$ is an stable interpolation set for $PW_K$. 
\end{itemize}
In particular, $\La_I$ is a universal set of stable sampling (resp. interpolation). 
\end{teo}

\subsection*{Paper guideline}  

Section \ref{preliminaries} is devoted to the preliminaries. With a brief review on the  LCA groups, the Schwartz-Bruhat space (a substitute to the Schwartz class  for LCA groups) and Riemann integrable sets, we set the notation and recall some technical material adapted to our approach. 
We also introduce Landau-Beurling densities in LCA groups following  \cite{neclca}.  As for the definitions of quasicrystals and CP-schemes  on LCA groups we rely on \cite{Meyer2}.

\medskip

In Section \ref{existence of lattice} we characterize all LCA groups containing complete CP-schemes and consequently, construct all possible quasicrystals.  
First, we deal with the easier case $G=\R^d\times \T^\ell$  which leads to Theorem \ref{existence sacps} proved in Subsection \ref{noimpo}. 
Then, we proceed to the proof of Theorem \ref{imposible} and, since it is long and technical we split it in several parts, all given in Subsection~\ref{impo}. 

\medskip

Finally, Section \ref{DSI} is devoted to the problem of universal sampling and interpolation. In Subsection \ref{Poissons}, we prove Theorem \ref{nl} and discuss the relation of quasicrystals with Poisson measures, and in 
Subsection  \ref{BeurlingQ} we show Beurling type Theorem \ref{thsuf} for simple quasicrystals.

\section{Preliminaries on LCA groups}\label{preliminaries}

Throughout this section we  review basic facts on  locally compact abelian (LCA) groups, setting in this way the notation we need for the following sections 
(for more details see \cite{DE},  \cite{HR1}, \cite{HR2}, and \cite{rudi}).  

\subsection{Basic definitions and notation}

Let $G$ denote a Hausdorff locally compact abelian  group, and $\widehat{G}$ its dual group, that is:
$$
\widehat{G}= \{\gamma: G \to \C, \ \mbox{and}\ \gamma\ \mbox{is a continuous character of}\ G\},
$$
where a character is a function satisfying the following properties:
\begin{itemize}
\item[(i)] $ |\gamma(x)| =1, \,  \forall x\in G$;
\item[(ii)] $\gamma(x+y)=\gamma(x)\gamma(y),  \, \forall x, y \in G$.
\end{itemize}

Thus, the characters generalize the exponential functions $\gamma(x)=\gamma_t(x)= e^{2\pi it x}$ in the case $G=(\R,+)$.  
Given $H$ a closed subgroup of $G$, the annihilator of $H$ is denoted by $H^\bot$, where
$$
H^\bot:= \{\gamma \in \G: \gamma(h)=1, \,\forall h\in H\}.
$$

\medskip

\begin{fed}
Let $G$ be a LCA group. A lattice $H$ of $G$ is a discrete subgroup such that $G/H$ is compact. 
\end{fed}

The annihilator $\Gamma=H^\bot$ is a lattice in $\G$ called \textbf{dual lattice}. On the other hand, it can be proved that there exists a relatively compact Borel set of representatives of $G/ H$. Any of these sets will be called  \textbf{section} of $G/H$ (see \cite{FGg} and \cite{KK}). 

\medskip

On every LCA group $G$ there exists a Haar measure; that is, a non-negative, Borel regular measure 
$\mu_G$ that is non-identically zero and translation-invariant, which means
$$
\mu_G(E+x)=\mu_G(E),
$$
for every element $x\in G$ and every Borel set $E\subset G$. This measure is unique up to a constant. Analogously to the Lebesgue spaces, we can define the $L^p(G)=L^p(G, \mu_G)$ spaces associated to the group $G$ and the measure $\mu_G$:
$$
L^p(G):= \Big\{f:G \to \C, \,\,\, f \peso{is measurable and} \int_G |f(x)|^p\,d\mu_G(x)<\infty\Big\}.
$$

\medskip

\begin{teo} Let $G$ be an LCA group and $\widehat{G}$ its dual. Then
\begin{itemize}
\item[(i)] The dual group $\widehat{G}$, with the operation $(\gamma+\gamma')(x)= \gamma(x)\gamma'(x)$ is an LCA group. The topology in $\widehat{G}$ is the one induced by the identification of the characters of the group with the characters of the algebra $L^1(G)$.
\item[(ii)] The dual group of $\widehat{G}$ is  isomorphic (as topological groups) to $G$, that is, 
$\widehat{\widehat{G \,}}\approx G$, with the identification $g\in G \leftrightarrow e_g \in \widehat{\widehat{G\,}}$, 
where $e_g(\gamma):=\gamma(g)$.
\item[(iii)] $G$ is discrete (resp. compact) if and only if $\widehat{G}$ is compact (resp. discrete).
\end{itemize}
\end{teo}

As a consequence of (ii) of the previous theorem, we could use the notation $(x,\gamma)$ for the complex number $\gamma(x)$, representing either the character $\gamma$ applied to $x$ or the character $x$ applied to $\gamma$.

\medskip

Taking $f\in L^1(G)$ we define the Fourier transform of $f$, as the function $\hat{f}:\widehat{G}\to\C$ given by
$$
\hat{f}(\gamma)=\int_G f(x)(x,-\gamma)\, d\mu_G(x), \,\,\, \gamma \in \widehat{G}.
$$
If the Haar measure of the dual group $\widehat{G}$ is normalized conveniently,  we obtain the inversion formula
$$
f(x)= \int_{\widehat{G}} \hat{f}(\gamma) (x,\gamma)d\mu_{\widehat{G}}(\gamma),
$$
for a specific class of functions. If the Haar measures $\mu_G$ and $\mu_{\widehat{G}}$ are normalized in such a way, the Fourier transform on $L^1(G)\cap L^2(G)$ can be extended to a unitary operator from $L^2(G)$ onto $L^2(\widehat{G})$. Thus the Plancherel formula holds:
$$
\pint{f,g}= \int_G f(x)\overline{g(x)} d\mu_G(x)
= 
\int_{\widehat{G}} \hat{f}(\gamma)\overline{\hat{g}(\gamma)} d\mu_{\widehat{G}}(\gamma) = \langle \hat{f},\hat{g}\rangle,
$$
for $f,g \in L^2(G)$.

\begin{fed} \label{mass}
Given a lattice $H$ of $G$, we denote by $\m{H}$ the total mass of a section with respect to some normalized Haar measure of $G$. 
\end{fed}

In what follows, we normalize the Haar measure so that the inversion formula holds. So, we get that $\m{H}= 1/\m{H^{\bot}}$.

\subsection{The Schwartz-Bruhat space $\swb$}\label{swb}

In \cite{Bru}, Bruhat extended the notion of $C^\infty$ function to a large class of locally compact groups, which includes the LCA groups. Inherent in this, there is an extension of the idea of the Schwartz class $\swb$ to any LCA group. The definition of $\swb$, as it was presented in \cite{Weil}, goes as follows. If $G $ is a Lie group, and so $G \simeq \R^d\times \Z^m\times \T^\ell\times F$,  with $F$ a finite abelian group, then a function $f$ is in $\swb$ if it lies in $C^\infty(G)$, and if $P(\partial)f$ remains bounded on $G$ for every polynomial differential operator $P(\partial)$. As usual, it can be proved that this is equivalent to requiring that $P(\partial)f\in L^2(G)$ for every polynomial differential operator $P(\partial)$. If G is any LCA group,
then
\begin{equation}\label{direct limit}
\swb=\varinjlim\,\ese(H/K),
\end{equation}
where the direct limit is over pairs $(H, K)$ of subgroups of $G$ such that $H$ is open and compactly generated, $K$ is compact, and $H/K$ is a Lie
group. Since $H$ is compactly generated, $H/K$ modulo its identity component is finitely generated; thus $H/K$ has the form $\R^d\times \Z^m\times \T^\ell\times F$. In \eqref{direct limit}, the set $\ese(H/K)$ should be understood as those continuous functions with support contained in $H$, that are constants in the coset modulus $K$, and their projections to the quotient $H/K$ belong to $\ese(H/K)$.

\medskip

An intrinsic characterization of $\swb$, without using direct limits, was proved by Osborne in \cite{O}. In order to present this characterization %of $\swb$,
we need the following definition.

\medskip

\begin{fed}\label{A-functions}
Define $\ag$ as the set of those $f\in L^\infty(G)$ for which there exists a compact set $C_f$ such that the following holds: 
for any positive integer $n$, there exists a constant $M_n>0$ such that for each integer $k\geq 1$,
$$
\|f\chi_{G\setminus C_f^{(k)}}\|_\infty\leq M_n k^{-n},
$$
where $C_f^{(1)}=C_f$ and $C_f^{(k)}=C_f^{(k-1)}+C_f$ for $k\geq 2$.
\end{fed}

\begin{rem}
Note that $C_f$ may be enlarged at will. Hence, we can assume that it is a neighborhood of the identity of $G$. \EOE
\end{rem}

In the following proposition we name some of the basic properties of $\ag$ (see ~\cite{O} for their proofs).
\begin{pro}
The class $\ag$ has the following properties:
\begin{enumerate}
\item[$i)$] It is translation invariant;
\item[$ii)$] It is closed under pointwise multiplication;
\item[$iii)$] It is a convolution algebra;
\item[$iv)$] lt is included in $L^1(G)$. 
\end{enumerate}
\end{pro}

Now we are ready to give the Osborne's characterization of the class $\swb$ (see ~\cite{O}).

\begin{teo}\label{Osborne identity}
The Schwartz-Bruhat class can be characterized in the following way:
 $$
 \swb=\{f\in\ag:\ \widehat{f}\in\agd\}.
 $$
\end{teo}

As a direct consequence of Osborne's characterization of the class $\swb$ we get:
\begin{cor}\label{tensor product} 
Let $G_1$ and $G_2$ two LCA groups, $f\in\swbp{G_1}$, and $g\in \swbp{G_2}$. Then, the function $h:G_1\times G_2\to\C$ defined by
$$
h(x_1,x_2)=f(x_1)g(x_2)
$$
belongs to $\swbp{G_1\times G_2}$.
\end{cor}

Recall that a trigonometric polynomial in the group setting is a finite linear combination of characters.  
Some useful properties of the Schwartz-Bruhat class are the following. 

\begin{lem}\label{vaaam}
Let $f\in \swbp{G}$ and let $p$ be a trigonometric polynomial. Then $f\cdot p\in\swbp{G}$. 
\end{lem}
\bdem
Since $p\in L^\infty(G)$, clearly $f\cdot p\in \ag$. On the other hand, $\widehat{f\cdot p}$ is a linear combination of translations of $\widehat{f}$, which belongs to $\agd$. Therefore,  $\widehat{f\cdot p}\in \agd$. 
\edem

\begin{pro}
Let  $H$ be a lattice of $G$ and let $\Gamma$ %=H^{\bot}
be the dual lattice. Then, given any function $f\in \swb$,  the Poisson's summation formula holds, which means that for every $x\in G$
\begin{equation}\label{eq Poisson}
\sum_{h\in H} f(x+h) = \frac{1}{\m{H}}\sum_{\gamma\in\Gamma} \widehat{f}(\gamma) e_x(\gamma), 
\end{equation}
where both series in \eqref{eq Poisson} are absolutely convergent, and the constant $\m{H}$ is given by Definition~\ref{mass}.
\end{pro}

\subsection{Riemann integrable sets}

In this section we recall some known properties of Riemann integrable sets to the setting of LCA groups. Since we could not find some of them in the literature, we will give their proofs.

\begin{fed}\label{riemann integrable set}
A measurable set $A\subseteq G$ is called Riemann integrable if the (Haar) measure of its (topological) boundary $\partial A$ is zero.
\end{fed}

Before going on, we point out that since the following technical results will be used 
for a compact set $\Omega\subseteq \G$, we will present them in this context.   

\begin{lem}\label{Urysohn}
Let $C$ be a compact subset of $\G$ and let $U$ be an open set of $\G$ with compact closure and such that $C\subseteq U$. Then, there exists a function $f\in\swbd$ such that
$$
\sub{\chi}{C} \leq f \leq  \sub{\chi}{U}.
$$
\end{lem}
\bdem
Let $V$ be a compact neighborhood of $e$ in $\G$ such that  $C+V\subseteq U$. Take a compact subgroup $K$ contained in $V$ which satisfies that
$$
\G/K\simeq \R^d\times \Z^m\times \T^\ell\times F,
$$
for some integers $d,m,\ell\geq 0$. Let $\pi:\widehat{G}\to\widehat{G}/K$ denote the canonical projection, and consider a compact neighborhood $\ewe$ of the identity in $\G/K$ such that $W=\pi^{-1}(\ewe)\subseteq V$. Choose a function $f_0\in\swbp{\G/K}$ such that 
$$
 \sub{\chi}{\pi(C)}\leq f_0\leq   \sub{\chi}{\pi(C+W)}.
$$
Then, clearly the function $f=f_0 \circ \pi$ belongs to $\swbd$ and it satisfies that
$$
 \sub{\chi}{C} \leq   \sub{\chi}{C+K} \leq f \leq    \sub{\chi}{C+W}\leq  \sub{\chi}{U},
$$
because  $\pi^{-1}(\pi(C))=C+K$, and $C+W=\pi^{-1}(\pi(C+W))$ since $W$ is saturated.
\edem

\medskip

Note that if $\Omega$ is Riemann integrable, then $L^p(\Omega)=L^p(\mbox{Int}(\Omega))$, where $\mbox{Int}(\Omega)$ denotes the interior of $\Omega$. Hence, as a consequence of Lemma \ref{Urysohn} and the regularity of the Haar measure $\mu_\G$, we obtain:
\begin{cor}
Let $\Omega$ be a relatively compact, Riemann integrable Borel subset of $\G$. Then
$$
\{f\in\swbd: \mbox{supp}(f)\subseteq \Omega\}
$$ 
is dense in $L^p(\Omega)$ for any $p\in [1,\infty)$.
\end{cor}

Another consequence of Lemma \ref{Urysohn} is the following characterization of Riemann integrable sets, which is well known in $\R^d$.
\begin{pro}\label{con Cinf} 
Let $\Omega\subseteq \G$ be a compact set. Then, the following statements are equivalent: 
\begin{enumerate}
\item The set $\Omega$ is Riemann integrable;
\item For every $\eps>0$, there exist positive functions $g_\eps, \,h_\eps\in \swbd$ such that 
$$
g_\eps\leq \sub{\chi}{\Omega}\leq h_\eps \peso{and} \int_{\G} (h_\eps-g_\eps)(\gamma)\,d\mu_{\G}(\gamma)\leq \eps.
$$
\end{enumerate}
\end{pro}
\bdem
$\ $
\begin{description}
\item[$1\Rightarrow 2)$] Let $U\subseteq \G$ be an open set with compact closure such that $\Omega\subseteq U$ and $\mu_\G(U\setminus \Omega)\leq \eps/2$. By Lemma \ref{Urysohn}, there exists a function $h_\eps\in\swbd$  such that 
$$
\sub{\chi}{\Omega} \leq h_\eps \leq \chi_U.
$$
On the other hand, if $\mu_\G(\Omega)=0$ then taking $g_\eps\equiv 0$ the implication is proved. If $\mu_\G(\Omega)\neq 0$, then $\mu_\G(\mbox{Int}(\Omega))\neq 0$. Let $\Omega_\eps$ be a compact subset of $\mbox{Int}(\Omega)$ such that $\mu_\G(\Omega \setminus \Omega_\eps)\leq\eps/2$, and let $g_\eps$ be a function in $\swbd$ that satisfies $\chi_{\Omega_\eps}\leq g_\eps\leq \chi_{\mbox{Int}(\Omega)}$. 
Then $g_\eps\leq \sub{\chi}{\Omega}\leq h_\eps$, and it holds that
$$
 \int_{\G} (h_\eps-g_\eps)(\gamma)\,d\mu_{\G}(\gamma)\leq \mu_\G(U\setminus \Omega_\eps)=
  \mu_\G(\Omega \setminus \Omega_\eps)+\mu_\G(U\setminus \Omega) \leq \eps.
$$
\item[$2\Rightarrow 1)$] It is enough to note that for every $\eps>0$, 
$$
 \mu_{\G}(\partial \Omega) \leq  \int_{\G} (h_\eps-g_\eps)(\gamma)\,d\mu_{\G}(\gamma)\leq \eps.
$$
\end{description} 
\edem

\begin{pro}\label{big brother}
Let $\Omega$ be a compact subset of $\G$ and $\eps>0$. Then, there exists a compact set $\Omega'\supseteq \Omega$ which is Riemann integrable and $\mu_{\G} (\Omega' \setminus \Omega)\leq \eps$.
\end{pro}
\bdem
Let $U\subseteq \G$ be an open set containing $\Omega$, and such that $\mu_{\G} (U \setminus \Omega)<\eps$. By Lemma \ref{Urysohn}, there is function $f\in\swbd$ such that 
$\chi_\Omega \leq f \leq \chi_U$. For each $\alpha\in(0,1)$ consider the level set $C_\alpha=\{\gamma \in \G : f(\gamma)=\alpha\}$. Since the sets $C_\alpha$ are pairwise disjoint, and all of them are contained in $U$, which has finite measure, at most a countable number of them can have positive measure. Take $\alpha_0$ such that $\mu_\G(C_{\alpha_0})=0$. Then, the set
$
\Omega'=\{\gamma \in \G :\ f(\gamma)\geq \alpha_0\}
$
satisfies the required properties.
\edem

\subsection{Landau-Beurling densities on LCA groups} \label{LBdensities}

We start the section with the definition of Landau-Beurling densities in LCA groups introduced in \cite{neclca}.   
In $\R^d$, these densities compare the concentration of the points of a given discrete set with that of  the integer lattice $\Z^d$. In a topological group, this comparison is done with respect to some reference lattice by means of the following relation (see \cite{neclca}). Recall that a subset $\La$ of $G$ is called \textbf{uniformly discrete} if there exists an open set $U$ such that the sets $\la+U$ ($\la$ in $\La$) are pairwise disjoints.

\begin{fed}
Given two uniformly discrete sets $\La$ and $\La'$ of an LCA group $G$,  and non-negative numbers $\alpha$ and $\alpha'$, we write $\alpha\La\leqm\alpha'\La'$ if for every $\eps>0$ there exists a compact subset $K$ of $G$ such that for every compact subset $L$ we have
$$
(1-\eps)\alpha\,\#(\La\cap L)\leq \alpha'\,\#(\La'\cap (K+L)).
$$ 
\end{fed}

Now, we have to fix a reference lattice in the group $G$.  We have assumed that $\widehat{G}$ is compactly generated, thus $G$ is isomorphic to $\R^d\times\T^m\times \D$, where $\D$ is a countable discrete group. So, a natural reference lattice is $H_0=\Z^d\times \{e\}\times \D$. Using this reference lattice, and the above transitive relation, we have all what we need to recall the definitions of upper and lower densities (see \cite{neclca} for further details).

\begin{fed}\label{densities LCA}
Let $\La$ be a uniformly discrete subset of $G$. The \textit{lower uniform density} of $\La$ is defined as
$$
\ede^-(\La)=\sup\{\alpha \in \R^+: \alpha H_0\leqm \La\}.
$$
On the other hand, its \textit{upper uniform density} is 
$$
\ede^+(\La)=\inf\{\alpha \in \R^+: \La\leqm\alpha H_0\}.
$$
\end{fed}

These densities always satisfy that $\ede^-(\La)\leq \ede^+(\La)$, and they are finite. Moreover, it can be shown that the infimum and the supremum are actually a minimum and a maximum. In the case that both densities coincide, we will simply write $\ede(\La)$. It should be also mentioned that  in the case of $\R^d$, these densities coincide with the Landau-Beurling's densities when the reference lattice is $\Z^d$.

\medskip

Using these densities, Gr\"ochenig, Kutyniok, and Seip obtained in \cite{neclca} the following extension of the classical result of Landau to LCA-groups.

\begin{teo}\label{LCA Landau}
Suppose $\La$ is a uniform discrete subset of $G$. Let $\mu_G$ be normalized such that any fundamental domain of $H_0$ has measure one, and $\mu_\G$ so that the inversion formula holds. Then
\begin{enumerate}
\item[S)]   If $\La$ is a stable sampling set for $PW_\Omega$, then $\cD^-(\La)\geq \mu_{\widehat{G}}(\Omega)$; 
\item[I)]    If $\La$ is an stable interpolation set for $PW_\Omega$, then $\cD^+(\La)\leq \mu_{\widehat{G}}(\Omega)$.
\end{enumerate}
\end{teo}

\subsection{Cut and project schemes and quasicrystals in LCA groups} 

Given a lattice $H$  in $\R^m\times G$ such that  canonical projections
$$
p_1: \R^m\times G \to \R^m \peso {and} p_2: \R^m\times G \to G
$$
satisfy the following conditions:
\begin{itemize}
\item[C1)]   $p_1(H) $ is dense in $\R^m$;  
\item[C2)]   $\left.p_2\right|_H$ is one to one.
\end{itemize}
We say that the triple  $(\R^m,G,H)$ is a \textbf{cut and project scheme} (CP-scheme). These schemes were introduced by Meyer in \cite{Meyer2}.
If, in addition, the scheme also satisfies 
\begin{itemize}
\item[C3)]   $\left.p_1\right|_H $ is one to one;  
\end{itemize}
it is called \textbf{aperiodic}. Finally, a \textbf{complete} CP-scheme is a scheme that sati-

\noindent sfies (C1), (C2), (C3) and
\begin{itemize}
\item[C4)]   $p_2(H) $ is dense in $G$.
\end{itemize}
Different from the Euclidean space,  the existence of complete (even aperiodic)  CP-schemes in an LCA group is not always guaranteed.  
The problem  resides in the existence of such a lattice $H$ and it is a non trivial issue. In fact, this will be the subject of the next section.

\medskip

Assume for the moment  that there is a complete CP-scheme $(\R^m,G,H)$, and let $S\subset\R^m$ be a Riemann integrable set with non-empty interior. Recall that, the quasicrystal associated to $(\R^m,G, H)$ and $S$, is given by 
\begin{equation*}
\label{(3.1)}
\L_S=\{p_2(h): h \in H, \,p_1(h)\in S\}, 
\end{equation*}
where  $p_1: H \to \R^m \,$ and  $\, p_2: H \to G$ are the canonical projections.

\medskip
 
Let  $\Gamma \subseteq \R^m \times \G $ be the dual lattice of $H$ and consider the canonical projections
$$
q_1: \R^m\times \G \to \R^m \peso {and} q_2: \R^m\times \G \to \G.
$$
It can be seen that 
the projections $ \left.q_1\right|_\Gamma \,$ and $ \,\left.q_2\right|_\Gamma$ are injective and have dense range on  
$\R^m$ and $\G$ respectively (see Lemma 2, pg. 41, \cite{Meyer2}). 

\medskip

Given a compact set  $K\subset \G$,  we define the set
\begin{equation*} %\label{dualm}
M_K := \{q_1(\gamma): \gamma\in \Gamma, \, q_2(\gamma) \in K \} \subseteq \R^m.
\end{equation*}

\medskip

Recall that a quasicrystal is called simple if $m=1$ and $S=I$ is an interval. 
Note that if $I=[-\a, \a],$ then the set $\L_I$ is symmetric, that is $\L_I=-\L_I$.
 
 \medskip

\begin{lem} \label{sep}
Let $S$ be a compact set of $\R^m$, and let $K$ be a compact set of $\G$. 
Then $\L_S$ as well as $M_{K}$ are uniformly discrete.
\end{lem}

\bdem 
We will prove that $\L_S$ is uniformly discrete. The proof for $M_K$ is essentially the same. Assume by 
contradiction that $\L_S$ is not uniformly discrete, and  let $d$ denote a translation invariant metric in $G$. 

\medskip

Let $n_1,m_1\geq 1$ so that
\begin{align*}
d_1&=d(x_{n_1},x_{m_1})\leq\frac{1}{2},  \peso{with}  x_{n_1},x_{m_1} \in \L_S.\\
\intertext{Inductively, we look for  $n_k,m_k\geq 0$ so that}
d_k&=d(x_{n_k},x_{m_k})\leq\frac{d_{k-1}}{2 },  \peso{with}  x_{n_k},x_{m_k} \in \L_S.
\end{align*}
Let $\mu_k\in H$ such that $p_2(\mu_k)=x_{n_k}-x_{m_k}$. Then, by the construction of $(x_{n_k},x_{m_k})$,
the elements  $\mu_k$ have the following properties:
\begin{itemize}
\item  $\mu_k \neq \mu_j$ if $j\neq k$;
\item $p_2(\mu_k)\xrightarrow[k\to\infty]{} e\in G$;
\item $p_1(\mu_k)\in S-S$.
\end{itemize}
Then, taking  $k_0$ big enough, we have that for $k>k_0$ the elements $\mu_k$  belong to $(\,C\times (S-S)\,) \cap H$, where $C$ is a compact neighborhood of the identity. This implies that the sequence $\{\mu_k\}$ has an accumulation point, which is a contradiction because  $H$ is a uniformly discrete set. Thus, $\La_S$ is uniformly discrete.
\edem

\section{Existence of quasicrystals } \label{existence of lattice} %CP-schemes
 
The main objective of this section is to prove Theorem \ref{imposible} stated in the Introduction. For those groups, we describe the structure  and construct essentially all possible lattices $H$ that produce a complete CP-scheme $(\R^m,G,H).$ 

\subsection{Initial comments} \label{initials}

Recall that we are working with an LCA group $G$, which is  isomorphic to $\R^d\times \T^\ell\times \D$, where $\D$ is a countable discrete group. The main obstacles in trying to construct a CP-scheme $(\R^m,G,H)$ are hidden in the structure of the discrete group $\D$. Indeed, as we will see in the next subsection, if $G$ is connected then it is not difficult to explicitly construct a complete CP-scheme $(\R^m, G,H)$. In particular, this will prove next  result.

\medskip

\begin{teo}\label{existence sacps}
If $m\in\N$ and $G_0=\R^d\times \T^\ell$,  then there exists a complete CP-scheme $(\R^m, G_0, H_0)$.
\end{teo}

Suppose now that $\D$ is not trivial and that $(\R^m,G,H)$ is an CP-scheme. If $H_j=p_j(H)$ for $j=1,2$, then
$\phi:H_2\to H_1$ defined by
$$
\phi(p_2(h))=p_1(h)
$$
is a group homomorphism. Therefore, the possible  existence of elements with finite order in $\D$ imposes a restriction.
To see this, assume that we have a lattice $H$ in $\R^m \times G$ where $G=\R^n\times\D$. Given $\si \in \D$, if we call $H_\si= \{(x,y,r)\in H: r=\si\}$, then
$H_0$ is a lattice in $\R^m\times\R^n\times\{0\}$ and $H_\si= H_0 + h_\si$ with $h_\si \in H_\si.$
So $H_\si$ is a translation of the lattice $H_0$ along the element $h_\si$.

\medskip

The strategy to construct an CP-scheme is to start with a lattice in $\R^m \times\R^n$ with the required properties,
and pick the right translations $h_\si$ for each of the elements in $\D.$ Now if $\si$ has finite order, say $s$,  in $\D$,  $h_\si$ has to have a special form.
For, note that if $h_\si\in H_\si$,  then $s. h_\si \in H_0$, that is $s.h_\si=(x,y,0) \in H_0,$
which implies that  $h_\si=(x/s,y/s,\si).$ In other words, the translation $h_\si$ corresponding to an element of order $s$,
has its first two components in the lattice $\frac{1}{s}\tilde{H}_0$ where $ \tilde{H}_0=\{(x,y)\in \R^m\times\R^d: (x,y,0) \in H\}.$

\medskip

This special structure produced by the finite order elements of $\D$ on the translations, imposes restrictions on the existence of CP-schemes.
Let us show this with a toy example. 

\begin{exa}
Assume that $G=\R\times \Z_2$, and take two real numbers $\alpha$, $\beta$ such that $1$, $\alpha$, and $\beta$ are rationally independent. We first choose in $\R \times \R$  the lattice 
$$
\tilde{H}_0=\big\{\big(n + \alpha m, n + \beta m \big): \ \mbox{$n,m\in\Z$}\big\}.
$$
Here, the special choice of $\alpha$ and $\beta$ to be rationally independent is to obtain the density of the projections,
i.e. properties C1, and C4. Now we consider the translations and construct the CP-scheme
$(\R,G,H)$,  where 
$$
H=\big\{\big((n+r/2)+\alpha m, (n+r/2)+ \beta m,r\big): \ \mbox{$n,m\in\Z$ and $r\in\{0,1\}$}\big\}.
$$
Analogously we can construct a lattice $H$ in $\R\times(\R\times(\Z_2)^2)$ in such a way that $(\R,\R\times\Z_2^2,H)$ is an CP-scheme. Indeed, consider the lattice $H$ given by
$$
H=\big\{\big((n+{r}/{2})+\alpha (m+s/2), (n+r/2)+ \beta (m+s/2),r,s\big): \mbox{$n,m\in\Z$  $r,s\in\{0,1\}$}\big\}.
$$
The problem appears  when we try to deal with a third copy of $\Z_2$. In that case, we have three $\Z-$independent elements of order two in $\D$, while we only have two copies of $\R.$ This will be incompatible with the injectivity of the projection $p_1,$ as it will be shown in the proof of Proposition \ref{si se pasa chau}.
\end{exa}

The study of this simple example leads us to the characterization of those $\D$ such that, given $m\in \N$, there exists an CP-scheme $(\R^m, \R^d\times\T^\ell\times \D, H)$, that is Theorem \ref{imposible}.  The proofs  of Theorems  \ref{existence sacps} and \ref{imposible} are quite technical we leave them for the next subsections.

\subsection{Proof of Theorem \ref{existence sacps}} \label{noimpo}

In order to prove Theorem \ref{existence sacps} we need the following two technical but simple lemmas and some auxiliary definitions. Given $\vec\alpha\in \R^m$, and $\vec\beta\in\R^d$ 
we define the matrix
\begin{equation} \label{laT}
T_{\vec\alpha,\vec\beta} :=\begin{pmatrix}
0_{m\times m}&\vec\alpha \cdot \vec\beta^{\,T}\\
\vec\beta\cdot \vec\alpha^{\,T} & 0_{d\times d}
\end{pmatrix},
\end{equation}
where the vectors in $\R^m$ and $\R^d$ are considered as column vectors, hence
$$
\vec\alpha \cdot \vec\beta^{\,T} =\begin{pmatrix}
\alpha_1\beta_1&\dots&\alpha_1\beta_d\\
\vdots&\ddots&\vdots\\
\alpha_m\beta_1&\dots&\alpha_m\beta_d
\end{pmatrix}. 
$$
Using that $\|\vec\alpha \cdot \vec\beta^{\,T}\|_{sp}=\|\vec{\alpha}\|_{\R^m}\cdot\|\vec{\beta}\|_{\R^d}$ we have:
\medskip

\begin{lem}
If the vectors $\vec\alpha$ and $\vec\beta$ satisfy that
$$
\|\vec{\alpha}\|_{\R^m}\cdot\|\vec{\beta}\|_{\R^d}<1,
$$
then the matrix $(I+T_{\vec\alpha,\vec\beta})$ is invertible.
\end{lem}
\bdem
Since  $T_{\vec\alpha,\vec\beta}$ is symmetric, we obtain:
\begin{align*}
\|T_{\alpha,\beta}\|_{sp}^2&=\left\|\begin{pmatrix}
0_{m\times m}&\vec\alpha \cdot \vec\beta^{\,T}\\
\vec\beta\cdot \vec\alpha^{\,T} & 0_{d\times d}
\end{pmatrix}\right\|_{sp}^2=\left\|\begin{pmatrix}
0_{m\times m}&\vec\alpha \cdot \vec\beta^{\,T}\\
\vec\beta\cdot \vec\alpha^{\,T} & 0_{d\times d}
\end{pmatrix}^2\right\|_{sp}\\
&=\left\|\begin{pmatrix}
\|\beta\|_{\R^d}^2(\vec\alpha\cdot\vec\alpha^T)&0\\
0&\|\alpha\|_{\R^m}^2 (\vec\beta\cdot\vec\beta^T)
\end{pmatrix}\right\|_{sp}\\
&\leq \|\vec{\alpha}\|_{\R^m}^2\cdot\|\vec{\beta}\|_{\R^d}^2\\&<1.
\end{align*}
Therefore $I+T_{\vec\alpha,\vec\beta}$ is invertible.
\edem

The next lemma will allow us to chose vectors with arbitrary small size, and whose coordinates will satisfy a convenient independence condition over the field of the rational numbers.

\begin{lem}\label{lem rationally independent}
Given $n\in\N$ and $\eps>0$, there exists a vector $\vec{\xi}\in \R^n$ such that $\|\vec{\xi}\|_{\R^n}<\eps$ and its coordinates satisfy that for any choice of 
$\sigma\in\{1,-1\}^n$ the numbers
\begin{equation}\label{la condition del lema}
1,\xi_1^{\sigma_1},\xi_2^{\sigma_2}\ldots,\xi_{n-1}^{\sigma_{n-1}},\xi_n^{\sigma_n}
\end{equation}
are linearly independent over the field of the rational numbers.
\end{lem}  
\bdem
Indeed,  if $\{p_n\}$ denote the sequence of prime numbers, we can take
$$
\vec{\xi}=\Big(\frac{1}{\sqrt{p_{s+1}}},\ldots,\frac{1}{\sqrt{p_{s+n}}}\Big),
$$
where $s\in\N$ is big enough so that $\|\vec{\xi}\|_{\R^n}<\eps$. The linear independence over $\Q$ of the sequence$\{\sqrt{p_n}\}$ is due to Besicovitch (see for instance \cite{Bes}). On the other hand, since $\sqrt{p_n}=p_n/\sqrt{p_n}$, clearly $\vec{\xi}$, also satisfies \eqref{la condition del lema}.
\edem

By the previous lemma, we choose $\vec\alpha\in\R^m$, $\vec\beta\in\R^d$, and $\vec\gamma\in\R^\ell$ so that for any choice of signs in the powers, the $(m+d+\ell)$ real numbers
\begin{equation}\label{la condition1}
1,\alpha_1^{\pm 1},\ldots,\alpha_m^{\pm 1}, \beta_1^{\pm 1},\ldots,\beta_d^{\pm 1},\gamma_1^{\pm 1},\ldots,\gamma_\ell^{\pm 1},
\end{equation} 
are $\Q$-linearly independent, and
\begin{equation}\label{la condition2}
\|\vec\alpha\|_{\R^m}\|\vec\beta\|_{\R^d}<1.
\end{equation}
Set
$
\vec{\omega}_\gamma=(e^{2\pi i \alpha_1\gamma_1},\ldots,e^{2\pi i \alpha_1\gamma_\ell}),
$
and for $a \in \N$ we denote 
$$
\vec{\omega}_\gamma^{\,a}=(e^{2\pi i \alpha_1\gamma_1a},\ldots,e^{2\pi i \alpha_1\gamma_\ell a}).
$$ 
Now we are ready  to prove Theorem \ref{existence sacps}.

\bdem[Proof of Theorem \ref{existence sacps}] 
We start defining the set 
\begin{equation}\label{eq el H0}
H_0:=\left\{\left((I+T_{\alpha,\beta})\begin{pmatrix}\vec{n}\\ \vec{k}\end{pmatrix},
\vec{\omega}_\gamma^{\,n_1}\right)\right\}_{\vec{n}=(n_1,\dots,n_m)\in\Z^m\,,\ \vec{k}=(k_1,\dots,k_d)\in\Z^d}.
\end{equation}
We will prove that $H_0$ is a lattice of $\R^m\times G_0$ such that $(\R^m, G_0, H_0)$ is a complete CP-scheme. By definition, clearly $H_0$ is a discrete subgroup of $\R^m\times G$. On the other hand, the quotient $(\R^m\times G_0)/H_0$ can be identified with
$$
(I+T_{\alpha,\beta})\big( [0,1)^{m+d}\big)\times \T^\ell 
$$
where in $(I+T_{\alpha,\beta})\big( [0,1)^{m+d}\big)$ we consider the standard addition in $\R^{m+d}$ modulus the lattice $(I+T_{\alpha,\beta})\Z^{m+d}$. So, the quotient $(\R^m\times G_0)/H_0$ can be identified with $ \T^{m+d+\ell}$ which is compact. This concludes the proof that $H_0$ is a lattice. Now, we will prove that $p_1$ is one to one on $H_0$. Assume that 
$$
p_1\left((I+T_{\alpha,\beta})\begin{pmatrix}\vec{n}\\ \vec{k}\end{pmatrix},\vec{\omega}_\gamma^{\,n_1}\right)=0
$$
for some $\vec{n}\in \Z^k$, and $\vec{k}\in\Z^d$ so that at least one of them is different from zero. Then we have that the system
\begin{align*}
n_1+\alpha_1\beta_1 k_1+\ldots+\alpha_1\beta_d k_d&=0\\
\vdots\hspace{2cm}&=\,  \vdots\\
n_m+\alpha_m\beta_1 k_1+\ldots+\alpha_m\beta_d k_d&=0,
\end{align*}
has a non trivial integer solution. Now multiplying both sides of the first equation by $\alpha_1^{-1}$, both sides of the second equation by $\alpha_2^{-1}$ and so on, we reach a contradiction with condition~\eqref{la condition1}. So, $p_1$ is one to one on $H_0$. In a similar way we can prove that $p_2$ is one to one. Indeed, note that it is enough to prove that the first coordinate of $p_2$ is one to one, which is exactly the same computation as before.  

In order to prove the density of $p_1(H_0)$ in $\R^m$, first note that $\Z^m\subset p_1(H_0)$. So, if $\pi_m:\R^m\to\R^m/\Z^m$ is the canonical projection, it is enough to prove that $\pi_m\circ p_1$ is dense in the quotient. For, if $\vec{e}_{1}$ denotes the vector of the canonical basis of  $\R^d$ with a one in the first coordinate, then when $x$ runs over the real numbers the $p_1$ projection of
$$
(I+T_{\alpha,\beta})\begin{pmatrix}\vec{0}\\ x\vec{e}_1\end{pmatrix}
$$
describes the line in the direction of the vector $\beta_1(\alpha_1,\ldots,\alpha_m)$. By condition \eqref{la condition1}, the numbers
$$
\beta_1^{-1}\,,\,\alpha_1\,,\,\ldots \,,\, \alpha_m
$$
are linearly independent over $\Q$. Therefore, the numbers
$$
1\,,\,\beta_1\alpha_1\,,\,\ldots \,,\, \beta_1 \alpha_m
$$
are also linearly independent over $\Q$. So, the image of 
$$
\left\{(I+T_{\alpha,\beta})\begin{pmatrix}\vec{0}\\ n\vec{e}_1\end{pmatrix}: n\in \Z \right\}
$$
by $\pi_m\circ p_1$ is dense in $\R^m/\Z^m$ by Kronecker's theorem. To prove the density of $p_2(H_0)$ in $\R^d\times\T^l$, 
we observe  that  $\Z^d\times\{1\}\subset p_2(H_0)$ where $\{1\}$ stands for the unity in $\T^l$. So, as before, it is enough to find a dense set in $(\R^d\times \T^\ell) / (\Z^d\times\{1\})\simeq \T^{d+l}$. Let $\vec{e}_1$ now denotes the vector of the canonical basis of  $\R^m$ with a one in the first coordinate. Then when $x$ runs over the real numbers $p_2$ projection of 
$$
(I+T_{\alpha,\beta})\begin{pmatrix} x\vec{e}_1\\ \vec{0} \end{pmatrix}
$$
describes in the $\R^d$ component the line in the direction of the vector 
$$
\alpha_1(\beta_1,\ldots,\beta_d).
$$
By condition \eqref{la condition1}, the numbers
$$
\alpha_1^{-1}\,,\,\beta_1\,,\,\ldots \,,\, \beta_d \,,\, \gamma_1\,,\,\ldots\,,\,\gamma_\ell
$$
are linearly independent over $\Q$. Therefore, the numbers
$$
1\,,\,\alpha_1\beta_1\,,\,\ldots \,,\, \alpha_1\beta_d\,,\,\alpha_1\gamma_1\,,\,\ldots\,,\,\alpha_1\gamma_\ell
$$
are also linearly independent over $\Q$. From now on, the computation is similar to the one done for $p_1$ using Kronecker's theorem. 
\edem

\subsection{Proof of Theorem  \ref{imposible}}  \label{impo}

The proof of Theorem \ref{imposible} is quite long. Therefore, we will divide it in several parts. Let us begin with the following technical lemma. It says basically that if we can construct an CP-scheme in $\R^m\times G$ with some properties, we can also construct a CP-scheme, with the same properties, when the discrete component of $G$ is replaced by a subgroup of $\D.$

\begin{lem}\label{quasiok con los subgrupos}
Let $m\in\N$,  $G=\R^d\times \T^\ell\times \D$, and $\widetilde{\D}$ a subgroup of $\D$. Suppose that $H$ is a lattice of $\R^m\times G$, and define
$$
\widetilde{H}=\{(x,g)\in H: g\in\widetilde{G}\},
$$
where $\widetilde{G}=\R^d\times \T^\ell\times \widetilde{\D}$. Then $\widetilde{H}$ is a lattice in $\R^m\times \widetilde{G}$. Moreover, if   $(\R^m,G,H)$  satisfies any subset of the conditions (C2), (C3), and (C4), then $(\R^m,\widetilde{G},\widetilde{H})$ also  satisfies the same subset of conditions. 
\end{lem}
\bdem
Clearly $\widetilde{H}$ is a discrete subgroup of $\widetilde{G}$. On the other hand, identifying
$$
\big((\R^m\times G)/\widetilde{H}\big)/\widetilde{\pi}(H) \peso{with} (\R^m\times G)/H,
$$
we have the following commutative diagrams
$$\xymatrix{
\R^m\times G \ar[r]^{\widetilde{\pi}}\ar[dd]_\pi  &(\R^m\times G)/\widetilde{H} \ar[ddl]_{\pi_0} \\
&\\       
(\R^m\times G)/H &    }
\xymatrix{
\R^m\times \widetilde{G} \ar[rr]^{\widetilde{\pi}}\ar[dd]_\pi &  &(\R^m\times \widetilde{G})/\widetilde{H} \ar[ddll]_{\pi_0} \\
                        &  & \\
\pi(\R^m\times \widetilde{G} ) &  &   }.
$$
In these diagrams $\pi$, $\widetilde{\pi}$, and $\pi_0$ denote the corresponding canonical projections. Note that, given $(\vec{x},\tilde{g})$ and $(\vec{y},\tilde{h})$ in $\R^m\times\widetilde{G}$, these elements are $H$-related if and only if they are $\widetilde{H}$-related. Thus, $\pi_0$ is a homeomorphism between  $(\R^m\times \widetilde{G})/\widetilde{H} $ and $\pi(\R^m\times \widetilde{G} )$. Note also that $\pi(\R^m\times\widetilde{G})$ is closed in the quotient, which is compact. So, $\pi(\R^m\times\widetilde{G})$ is compact, and the same holds for $(\R^m\times\widetilde{G})/\widetilde{H}$. This shows that $\widetilde{H}$ is a lattice of $\R^m\times \widetilde{G}$. Clearly, if $p_1$ and/or $p_2$ restricted to $H$ is injective, the same happens with the corresponding restrictions to $\widetilde{H}$. On the other hand, since $\D$ is discrete, if $p_2(H)$ is dense, then for every $\si \in\D$
$$
p_2(H)\cap (\R^d\times \T^\ell\times \{\si\})
$$
is dense in the fiber $ \R^d\times \T^\ell\times \{\si\}$. Therefore, $p_2(\widetilde{\D})$ is also dense in $\widetilde{G}$. This completes the proof.
\edem

A direct consequence of this lemma is the following result:

\begin{cor}\label{se traslada}
Let $m\in\N$,  $G=\R^d\times \T^\ell\times \D$, and suppose that $H$ is a lattice in $\R^m\times G$. For each $r\in \D$, let $H_r$ be the subset of $H$ whose elements have coordinate in $\D$ equal to $r$. Then:
\begin{enumerate}
\item[i.)] $H_0$ is a lattice;
\item[ii.)] $H_r$ is a translation of $H_0$. 
\end{enumerate}
\end{cor}
\bdem
By Lemma \ref{quasiok con los subgrupos}, $H_{0}$ is a lattice in $\R^d\times \T^\ell\times \{0\}$. On the other hand, if $h_1,h_2\in  H_r$, then $h_1-h_2\in H_0$. 
\edem

As a particular case of these results, we get the following characterization of the lattices in $\R^n\times \T^\ell$ which is probably known.

\begin{lem}\label{lattices en rxt}
Given $n,\ell\in \N_0$, with at least one of them different from zero, then a lattice in $\R^n\times \T^\ell$ has the form
$$
\Big\{ \big(A\vec{k},\{e^{2\pi i (\,\vec{\gamma}_j\,\cdot\, \vec{k})}\,\}_{j=1}^\ell\,\big): \vec{k}\in \Z^n  \Big\}
$$
for some invertible $n\times n$ real matrix $A$ and vectors $\vec{\gamma}_1,\ldots,\vec{\gamma}_\ell\in\R^n$.
\end{lem}
\bdem[Sketch of the proof.] It is easier to describe the lattices in the dual space $\R^n\times \Z^\ell$. Let $E$ be one of these lattices, and given $\vec{k}\in\Z^\ell$, denote by $\La_{\vec{k}}=\{\la\in \R^m: (\la,\vec{k})\in E\}$ 
and denote by  $E_{\vec{k}} = \La_{\vec{k}} \times \{\vec{k}\}$ the subset of $E$ whose $\Z^\ell$-coordinate is precisely $\vec{k}$. By Corollary \ref{se traslada}, $E_{\vec{0}}$ is a lattice in $\R^n\times\{0\}$ and then  $\La_{\vec{0}}$ is a lattice in $\R^n$. So, there exists an invertible matrix $B$ such that
$$
\La_{\vec{0}}=B \Z^n.
$$
    On the other hand, if $(\la_1,\vec{k}),(\la_2,\vec{k})\in E_{\vec{k}}$, then $\la_1-\la_2 \in \La_{\vec{0}}$. Taking representatives $\la_{\vec{k}}\in \La_{\vec{k}}$ we have that $E_{\vec{k}}=(\la_{\vec{k}},\vec{k})+E_0$ and then we conclude that $E$ has the following structure
$$
E=\Big\{(B\vec{m}+\la_{\vec{k}}\,,\,\vec{k}):\ \vec{m}\in\Z^n\ \ \vec{k}\in\Z^\ell \Big\}.
$$
Considering this structure of $E$ and taking its dual we get the desired result.
\edem

Note that, in the lattice considered in Theorem \ref{existence sacps}, if $\vec{e}_1$ denotes the element of $\R^n$ that has one in the first coordinate and zero in the others, then
$$
A=I+T_{\alpha,\beta}\peso{and} \vec{\gamma}_j=\alpha_1\gamma_j\vec{e}_1.
$$

\medskip

The following proposition, combined with Lemma \ref{quasiok con los subgrupos}, proves one of the implications of Theorem \ref{imposible}.

\begin{pro}\label{si se pasa chau}
Let $m\in\N$ and let $G=\R^d\times \T^\ell\times \Z_p^{m+d+1}$. 
Then, it doesn't exist a lattice $H$ in $\R^m\times G$ such that $(\R^m, G, H)$ is a complete CP-scheme. 
%In particular, there is no CP-scheme of the form $(\R^m,G,H)$.
\end{pro}
\bdem
Suppose that there exists a complete CP-scheme $(\R^m,G,H)$. Let $\vec{u}_j\in  \Z_p^{m+d+1}$ be the element that is one in the $j$-th coordinate and zero in the other ones, and define
\begin{align*}
H_j&=\{(\vec{x},\vec{y},\vec{\w},\vec{r})\in H: \vec{r}=\vec{u}_j\},\\
\intertext{and}
H_0&=\{(\vec{x},\vec{y},\vec{\w},\vec{r})\in H: \vec{r}=\vec{0}\}.
\end{align*}
By condition (C4), for every $j$ it holds that $H_j\neq \varnothing$. So, take for every $j$
$$
v_j:=(\vec{x}_j,\vec{y}_j,\vec{\w}_j,\vec{u}_j)\in H_j.
$$
Then, for every $j$
$$
pv_j=(p\vec{x}_j,p\vec{y}_j,\vec{\w}_j^p,0)\in H_0.
$$
(we recall here the notation $\vec{\w}_j^p=(\vec{\w}_{j1}^p,\dots\vec{\w}_{j,\ell}^p) \in \T^{\ell}).$
By Lemma \ref{quasiok con los subgrupos}, $H_0$ is a lattice of $\R^m\times\R^d\times \T^\ell$ that satisfies the conditions $(C_3)$ and $(C_4)$. 
On the other hand, by Lemma \ref{lattices en rxt}, $H_0$ has a special form which implies that the vectors
$$
p\, \vec{x}_1,\ldots,p\, \vec{x}_{m+k+1}
$$
span a rational subspace of $\R^m$ of (rational) dimension at most $m+d$. This implies that there exist integers $n_j$, with no common divisors,  such that
$$
n_1\vec{x}_1+\ldots +n_{m+d+1}\vec{x}_{m+d+1}=0.
$$
Since the numbers $n_j$ do not have a common factor, there have to be one of them which is not multiple of $p$. Then
$$
v=n_1v_1+\ldots +n_{m+d+1}v_{m+d+1}\in H\setminus\{0\}
$$
and $p_1(v)=0$. This contradicts condition (C3), that is, the injectivity of $\left.p_1\right|_H$.
\edem

In order to prove the converse in Theorem \ref{imposible}, we use the following fact. Let $\D$ be a countable abelian group. If $\D$ does not have a copy of $\Z_p^{m+d+1}$ for any prime $p$, then it is isomorphic to a subgroup of a group of the form
\begin{equation}\label{eq los divisibles version corta}
\Q^{k_q}  \oplus   \bigoplus_{p: \text{ prime} }\Z(p^{\infty})^{k_p},
\end{equation}
where the exponents $k_p$ and $k_q$ are non-negative integers, $k_p\leq m+d$, and $\Z( p^{\infty} )$ are the so called quasi-cyclic groups defined by 
$$
\Z( p^{\infty} ) := \left\{r\in [0,1): r=\frac{k}{p^n}; \, k\in \N_0\,,\ n\in \N \right\}.
$$
We leave the proof of this fact to the Appendix \ref{App} (Proposition \ref{como son}). Now, motivated by this alternative description of countable abelian groups without copies of $\Z_p^{m+d+1}$ we consider the following result.

\begin{teo}\label{lattice en varios divisibles}
Let $m\in\N$, and let $G=\R^d\times \T^\ell\times \D$ where $\D$ is a divisible group of the form
$$
\D\simeq \Q^{k_q}  \oplus   \bigoplus_{p: \text{ prime} }\Z(p^{\infty})^{m+d}.
$$
Then, it is possible to construct a lattice $H$ in $\R^m\times G$ such that the CP-scheme $(\R^m, G, H)$ is complete. 
\end{teo}

\bdem 

We will start defining the lattice $H$ in $\R^m\times G$ that produce the complete CP-scheme. 
For this consider first 
$$
H_0=\left\{\left((I+T_{\alpha,\beta})\begin{pmatrix}\vec{n}\\ \vec{k}\end{pmatrix},\vec{\omega}_\gamma^{\,n_1}\right): \vec{n}\in\Z^m\,,\ \vec{k}\in\Z^d\right\}
$$
the lattice constructed in the proof of Theorem \ref{existence sacps}. 
Basically, $H$ will consists of some particular translations of the lattice $H_0$ located in each fiber
$$
\R^d\times\T^\ell\times\{r\}.
$$

Given $r\in\D$, write  $r = (r_\Q,r_P)$, where $r_\Q\in \Q^{k_q}, \; r_Q=\{r_Q(j):j=1,\dots,k_q\}$ and 
$$
r_P=\{r_{p}\}\in \bigoplus_{p: \text{ prime} }\Z(p^{\infty})^{m+d}.
$$ 
For each prime $p$,  $\;r_p\in \Z(p^{\infty})^{m+d}, \; r_p=\{r_p(1),\dots,r_p(m+d)\}, \; \text{ with } r_p(j) = {k_j}/{p^{s_j}}$ . 

With this conventions  define for $r\in\D$,  $ \vec{n}\in\Z^m$ and $\vec{k}\in\Z^d$,
\begin{align*}
\vec{u}(\vec{n},r)&=\left(\left(n_1+\sum_{p: \text{ prime} } r_p(1)\right),\ldots,\left(n_m+\sum_{p: \text{ prime} } r_p(m)\right)\right)\,,\\
\intertext{and}
\vec{v}(\vec{k},r)&=\left(\left(k_1+\sum_{p: \text{ prime} } r_p(m+1)\right),\ldots,\left(k_d+\sum_{p: \text{ prime} } r_p(m+d)\right)\right).
\end{align*}
Using these notations, we define $H$ as the set consisting of the elements of the form
\begin{equation}\label{lalatticeH}
\left((I+T_{\alpha,\beta})\begin{pmatrix}\vec{u}(\vec{n},r)\\ \vec{v}(\vec{k},r)\end{pmatrix}+\left(\sum_j r_\Q(j)\,\alpha_1\,\eta_j\right)\begin{pmatrix}\vec{e}_1\\\vec{0}\end{pmatrix}, \,\vec{\omega}_\gamma^{\,\vec{u}(\vec{n},r)_1},\,r\right)
\end{equation}
where $ \vec{n}\in\Z^m$, $\vec{k}\in\Z^d$, $r=(r_\Q,r_P)\in\D$,  $\vec{e}_1$ denotes the first element of the canonical basis of $\R^m$, $r\in \D$, 
and $\eta_1,\dots,\eta_{k_q}$ are real numbers which are rationally independent with $1$ and all the $\alpha$-s, $\beta$-s, $\gamma$-s and their inverses. The existence of the $\alpha$-s, $\beta$-s, $\gamma$-s and $\eta$-s is guaranteed by Lemma \ref{lem rationally independent}.

Straightforward computations show that $H$ is a lattice. Since $(\R^m,(\R^d\times \T^\ell), H_0)$ is a complete CP-scheme, clearly $p_1(H)$ and $p_2(H)$ are dense. 

In order to prove that the CP-scheme $(\R^m, G, H)$ is complete we only need to see that $p_1$ and $p_2$ restricted to  $H$ are one to one. 
For $p_1$  we proceed as in the proof of Theorem \ref{existence sacps}. Suppose that 
$$
p_1\left((I+T_{\alpha,\beta})\begin{pmatrix}\vec{u}(\vec{n},r)\\ \vec{v}(\vec{k},r)\end{pmatrix}+\left(\sum_j r_\Q(j)\,\alpha_1\,\eta_j\right)\begin{pmatrix}\vec{e}_1\\\vec{0}\end{pmatrix},\vec{\omega}_\gamma^{\,\vec{u}(\vec{n},r)_1},r\right)=0.
$$
Then
\begin{align*}
\vec{u}(\vec{n},r)_1+\alpha_1\beta_1 \vec{v}(\vec{k},r)_1+\ldots+\alpha_1\beta_d \vec{v}(\vec{k},r)_d&=-\sum_j r_\Q(j)\,\alpha_1\,\eta_j\\
\vec{u}(\vec{n},r)_2+\alpha_2\beta_1 \vec{v}(\vec{k},r)_1+\ldots+\alpha_2\beta_d \vec{v}(\vec{k},r)_d&=0\\
\vdots\hspace{2cm}&=\,  \vdots\\
\vec{u}(\vec{n},r)_m+\alpha_m\beta_1 \vec{v}(\vec{k},r)_1+\ldots+\alpha_m\beta_d \vec{v}(\vec{k},r)_d&=0.
\end{align*}

\noi The hypothesis of rational independence over the numbers $\alpha$-s, $\beta$-s, $\gamma$-s, and $\eta$-s implies, on the one hand that $r_\Q(j)=0$ for every $j$. On the other hand, it also implies that
\begin{align*}
\vec{u}(\vec{n},r)_1= \vec{v}(\vec{k},r)_1=\ldots=\vec{v}(\vec{k},r)_d&=0\\
\vdots\hspace{2cm}&=\,  \vdots\\
\vec{u}(\vec{n},r)_m=\vec{v}(\vec{k},r)_1=\ldots=\vec{v}(\vec{k},r)_d&=0.
\end{align*}

\noi Clearly this implies that the components in $\R^d$ and $\T^\ell$ are zero and one respectively. It only remains to prove that $r_P=0$. 

Assume by contradiction that $r_P\neq 0$.
Let $r_{\rho_1},\dots,r_{\rho_h}$ the non-zero components of $r_P$ where  $r_{\rho_t} \in \Z(\rho_t^{\infty})^{m+d}$  
for each $t=1,\dots,h$ and  the primes $\rho_t$ are all different.

Then, there exist non-negative integers  $k_{t,j} \text{ and positive integers } s_{t,j},$ with  $t = 1,\dots,h\;$  and  $j =1,\dots, m+d$, such that
$$r_{\rho_t}(j) = \frac{k_{t,j}}{\rho_t^{s_{t,j}}}.$$
The integers $k_{t,j} $ can be choosen in the range $0\leq k_{t,j} < \rho_t^{s_{t,j}} .$

Fix now $j \in \{1,\dots,m\},$ and set $k_{t} = k_{t,j}$ and $s_{t} = s_{t,j}.$ Since $\vec{u}(\vec{n},r)_{j}=0$, it holds that
$$
-n_j=\sum_{p: \text{ prime} } {r_p(j)} = \sum_{t=1}^h r_{\rho_t}(j) = \sum_{t=1}^h \frac{k_{t}}{\rho_t^{s_{t}}}.
$$
Then we have,
\begin{equation}\label{producto}
n_{j} \; \prod_{t=1}^h  \rho_t^{s_t} + \sum_{t=1}^h k_{t} \prod_{r\neq t} \rho_{r}^{s_{r}}=0.
\end{equation}

Assume now that one of the $k_{t}$ is different from zero, say $k_{1}\neq 0.$
Then rearranging terms in \eqref{producto} we have:
$$
- k_{1} \prod_{r\neq 1} \rho_{r}^{s_{r}} = n_{j} \; \prod_{t=1}^h  \rho_t^{s_t} + \sum_{t= 2}^h k_{t} \prod_{r\neq t} \rho_{r}^{s_{r}} =
\rho_1^{s_1}\left( n_{j} \; \prod_{t=2}^h  \rho_t^{s_t} + \sum_{t= 2}^h k_{t} \prod_{r\neq 1,t} \rho_{r}^{s_{r}} \right).
$$
So $\rho_1^{s_1}$ is a factor of $k_{1}$ which is a contradiction since $0\leq k_{1} < \rho_1^{s_{1}}$ and was assumed not to be zero.
Consequently, $r_{\rho_t} = 0$ for all $t$. This is true for every $j\in \{1,\dots,m\}$. A similar argument shows that the same holds for $j \in \{m+1,\dots m+d\}.$
We conclude that $r_P=0$ and then the projection $p_1$ is one to one on $H.$
\medskip

\noi Finally, if 
$$
p_2\left((I+T_{\alpha,\beta})\begin{pmatrix}\vec{u}(\vec{n},r)\\ \vec{v}(\vec{k},r)\end{pmatrix}+\left(\sum_j r_\Q(j)\,\alpha_1\,\eta_j\right)\begin{pmatrix}\vec{e}_1\\\vec{0}\end{pmatrix},\vec{\omega}_\gamma^{\,\vec{u}(\vec{n},r)_1},r\right)=0,
$$
then $r=0$. So the element of the lattice belongs to $H_0$. But there, $p_2$ is one to one. This concludes the proof.
\edem

\medskip

\noi As a consequence of Theorem \ref{lattice en varios divisibles} and Lemma \ref{quasiok con los subgrupos} we get the converse implication.

\begin{teo}\label{la voltereta} 
Let $m\in\N$ and let $G=\R^d\times \T^\ell\times \D$. If $\D$ does not have a copy of $\Z_p^{m+d+1}$ for any prime $p$, then there exists a complete CP-scheme $(\R^m, G, H)$.
\end{teo} 
\bdem
If $\D$ does not have a copy of $\Z_p^{m+d+1}$ for any prime $p$, then it is homeomorphic to a subgroup of a divisible group $\D^*$ of the form \eqref{eq los divisibles version corta} such that $k_p\leq m+d$ for every prime $p$. Let $G^*=\R^d\times \T^\ell\times \D^*$. By Theorem \ref{lattice en varios divisibles}, there exists a complete CP-scheme $(\R^m,{G^*},{H^*})$. By Lemma \ref{quasiok con los subgrupos}, this scheme induces a CP-scheme 
$(\R^m,\widetilde{G},\widetilde{H})$ that satisfies conditions (C2), (C3), (C4). Since it clearly satisfies (C1) too, it turns out that $(\R^m,\widetilde{G},\widetilde{H})$ is an CP-scheme.
\edem

\section{Universal stable sampling and interpolation}\label{DSI} 

\noindent 
This section is devoted to the construction of universal sets of stable sampling and interpolation in LCA groups, which we assume isomorphic to 
$\R^d\times \T^\ell\times \D$. 

\medskip

Suppose for a moment that $\ell=0$. Since we already know how to construct universal sets in $\R^d$, a
natural idea to construct such sets  in $G=\R^d\times \D$ could be to choose a universal set in each component.
The problem is how to combine them in order to get a universal set of stable sampling (resp. interpolation) in $G$.
For instance, consider a universal set of stable sampling $\La_0$ for $\R^d$ and then take in each connected component of $G$ a copy of $\L_0$. The resulting set will be a set of stable sampling for all $PW_\W$ such that the measures of almost all sections
$$
\Omega_{t_0}=\{(x,t)\in\W: t=t_0\in \widehat{\D}\}
$$ 
are bounded by the density of $\La_0$. As the measure of these sections can be much bigger than the measure of $\Omega$, this approach fails. Hence, in order to get universal sets of stable sampling or stable interpolation we need to combine universal sets in each components in a more delicate way. We are able to solve the problem using quasicrystals. The main weak point of this approach, is that quasicrystals do not always exist. However, we already know that there are many groups where we can find quasicrystals. 
In these groups, we show that quasicrystals satisfy universality properties. We essentially follow the ideas in \cite{MM1} in order to construct universal stable sampling and stable interpolation sets, adapted to the group setting.  The main difference in our context is the proof of Theorem \ref{nl}. Once this theorem is established, the duality theorem (Theorem \ref{duality}) can be proved in the same way as in the case of $\R^d$. So, we will give the complete proof of Theorem \ref{nl}, and refer the interested reader to \cite{MM1} for the details of the proof of Theorem \ref{duality}. The end of the section is devoted to the Beurling type theorem for simple quasicrystals. From this we will deduce the existence of universal sets of stable sampling and stable interpolation in the groups that admit simple quasicrystals.

\subsection{Poisson measures associated to quasicrystals}  \label{Poissons}

The main goal of this section is to provide the proof of Theorem \ref{nl}. 
Before proceeding to the proof, let us make some comments.

\subsection*{Density of $M_K$} \label{relations}

First of all, note that using Theorem \ref{nl} and Proposition \ref{con Cinf} we directly get the following result.

\begin{cor}\label{nl riemann}
Let $H$  be a lattice on $ \R^m \times G$ such that $\left.p_1\right|_H$ and $\left.p_2\right|_H$ are one to one, and $p_1(H)$ is dense in $\R^m$. If $J$ and $K$ Riemann integrable subsets of $\R^m$  and $G$ respectively, then
$$
\lim_{r\to \infty}\,  \frac{1}{r^m} \,  \sum_{h \in H} \chi_J\left(\frac{p_1(h)-a}{r}\right)\chi_K(p_2(h))
\, = \, \frac{m_{\R^m}(J) \ m_G(K)}{\m{H}} 
$$
uniformly in $a\in \R^m$. 
\end{cor}

As a consequence of this corollary we obtain the following result on the 
Beurling density of $M_K$.

\begin{cor}\label{uno}
Given a Riemann integrable compact set $K\subseteq \G$, it holds that
$$
\cD(M_K) = \frac{\mu_{\G}(K)}{\m{\Gamma}},
$$
where $\Gamma$ is the lattice related to $M_K$ by the cut and project method.
\end{cor}

\bdem
Using Corollary \ref{nl riemann}  with $J={[0,1]}^m$, 
$q_1,q_2$ instead of $p_1,p_2$, and $\Gamma$ instead of $H$, we get that
$$
\lim_{r\to \infty}\,  \frac{1}{r^m} \,\sum_{\gamma\in\Gamma}\sub{\chi}{J}
\left(\frac{q_1(\gamma)-a}{r}\right)\sub{\chi}{K}(q_2(\gamma))= \frac{1}{\m{\Gamma}} \cdot \mu_{\G}(K).
$$
 Since the limit is uniform in $a$, this precisely says that 
$\cD(M_K) = \m{\Gamma}^{-1}\mu_{\G}(K)$.  
\edem

\subsection*{Theorem \ref{nl} and Poisson measures}\label{PM}

Theorem \ref{nl} and Corollary \ref{nl riemann} are related with Poisson measures.  Let $H$ be a lattice in $\R^m\times G$ such that $\left.p_1\right|_H$ y $\left.p_2\right|_H$ are one to one and their images are dense 
in $\R^m$ and  $G$ respectively. Fix a function $\psi\in \swb$ whose Fourier transform is supported in the Riemann integrable compact set $K \subset \G$, and consider the atomic measure in $\R^m$ defined by
$$
\nu_\psi=\sum_{h\in H} \psi(p_2(h)) \delta_{p_1(h)}.
$$
The distributional Fourier transform of $\nu$ is another atomic measure on $\R$, whose formula is 
$$
\widehat{\nu}_\psi= \frac{1}{\m{H}}\sum_{\gamma\in \Gamma} \widehat{\psi}(q_2(\gamma)) \delta_{-q_1(\gamma)},
$$
where $\Gamma$ is the dual lattice of $H$, and $q_j$ are the dual projections of $p_j$ for $j=1,2$.  Note that $\widehat{\nu}$ is supported in the quasicrystal $M_K$. By Theorem \ref{nl}, both $\nu$ and $\widehat{\nu}$ are translation bounded measures, i.e.
$$
\sup_{x\in\R^m} \nu_\psi([x,x+1]^m)<\infty \peso{and} \sup_{x\in\R^m} \widehat{\nu}_\psi([x,x+1]^m)<\infty,
$$
and 
$$
\eme(\nu_\psi):=\lim_{r\to\infty} \frac{1}{r} \nu_\psi([-r,r])= \frac{1}{ \m{H} }\int_G \psi(g)\,dm_G.
$$
Therefore, they are Poisson measures in the sense of \cite{Meyer6}. In particular, both of them are almost periodic measures. 
Atomic measures whose (distributional) Fourier transform is also an atomic measure have been widely studied in connection with the so called Poisson type formulae (see \cite{Cordoba}, \cite{Lag}, \cite{K} and \cite{LevOlev}). 

\subsection*{Proof of Theorem \ref{nl}}

We will divide the proof of Theorem \ref{nl}  in two lemmas. Recall that  $\G$ is a compactly generated group and let 
$\Gamma\subset  \R^m \times \G$ be a uniform lattice. Since $\G$ is compactly generated we have that 
\begin{equation}\label{cg}
\G  \simeq \R^d  \times   \Z^{\ell}   \times   \K.
\end{equation}
Consider the lattice 
\begin{equation}\label{Ksi}
\Xi= \Z^d \times \Z^{\ell} \times \{e_{\K}\},
\end{equation}
where $e_{\K}$ denotes the identity of $\K$. Let  $Q=[0,1)^d \times \{0\}^\ell \times \K$, and for $\tau\in\G$ define
\begin{equation} \label{Qtau}
Q_{\tau}= Q+\tau. %, \,\,\, \xi \in \Xi.
\end{equation}
Note that $\G = \bigcup_{\xi \in \Xi} Q_{\xi}.$ With these notations we have the following results.

\begin{lem}\label{separacion uniforme}
Let $\G$ be a compactly generated group. For any $\tau \in \G$, the quasicrystal
$$
M_{Q_{\tau}}=\{q_1(\gamma):\gamma \in \Gamma,\, q_2(\gamma) \in Q_\tau\},
$$ 
is uniformly discrete and the separation constants do not depend on $\tau$.
\end{lem}

\bdem
Let $S\subseteq \G$ be any compact set such that $Q$ is contained in the interior of $S$. By Lemma \ref{sep}, the quasicrystal $M_{(S-S)}$ is uniformly discrete, where as usual
$$
S-S=\{s_1-s_2:\ s_1,s_2\in S\}.
$$ 
Let $\delta>0$ such that $|x-y|>\delta$ if $x,y\in M_{(S-S)} \subset \R^m$. Since, $e\in S-S$, 
note that $|x|>\delta$ for every non zero element of  $M_{(S-S)}$.

\medskip

Although it is not necessary, for the sake of clarity, let us firstly assume that $\tau \in q_2(\Gamma)$.  So, $\tau= q_2(\gamma)$, for some $\gamma \in \Gamma$. Given two different elements $x_1,x_2  \in M_{Q_\tau}$, there exist $\gamma_1,\gamma_2$ such that  $q_1(\gamma_i)=x_i$, $i=1,2$. So, if we define $\eta_i=\gamma_i -\gamma$  then  $q_2(\eta_i) \in Q$. Therefore
 $$
 |x_1-x_2|=|q_1(\gamma_1)-q_1(\gamma_2)| =|q_1(\eta_1-\eta_2)| \geq \delta,
 $$
because $q_2(\eta_1-\eta_2)\in Q-Q\subseteq S-S$.
 
\medskip
 
Now, suppose that $\tau \notin q_2(\Gamma)$, and take two different elements $x_1,x_2  \in M_{Q_\tau}$. Since $\left.p_1\right|_H$ is one to one, by a simple duality argument, we can see that $q_2(\Gamma)$ is dense in $\G$. So, there exists $\gamma\in\Gamma$ ``closed'' enough to $\tau$, such that $x_1,x_2\in S+q_2(\gamma)$, Thus, a similar argument as before shows that $x_1-x_2=q_1(\eta)$ for some $\eta\in \Gamma$ satisfying that $q_2(\eta)\in S-S$. Hence, $x_1-x_2\in M_{(S-S)}$, and since it is non zero
$$
 |x_1-x_2|\geq \delta,
 $$
which concludes the proof.
\edem

\begin{lem} \label{limit}
Let $\hat{\fii}\in\ese(\R^m) $ and $\hat{\psi}\in\ese(\G)$. Then
\begin{align}
\sum_{  \gamma \in\Gamma: \gamma \neq e}|\hat{\fii}\big(r\,q_1(\gamma)\big)\hat{\psi}\big(q_2(\gamma)\big)|\xrightarrow[r\to\infty]{}0\label{que pollo},
\end{align}
where $e$ denotes the identity of $\G$.
\end{lem}

\begin{figure}[H]
           \centering
           \includegraphics[height=6cm]{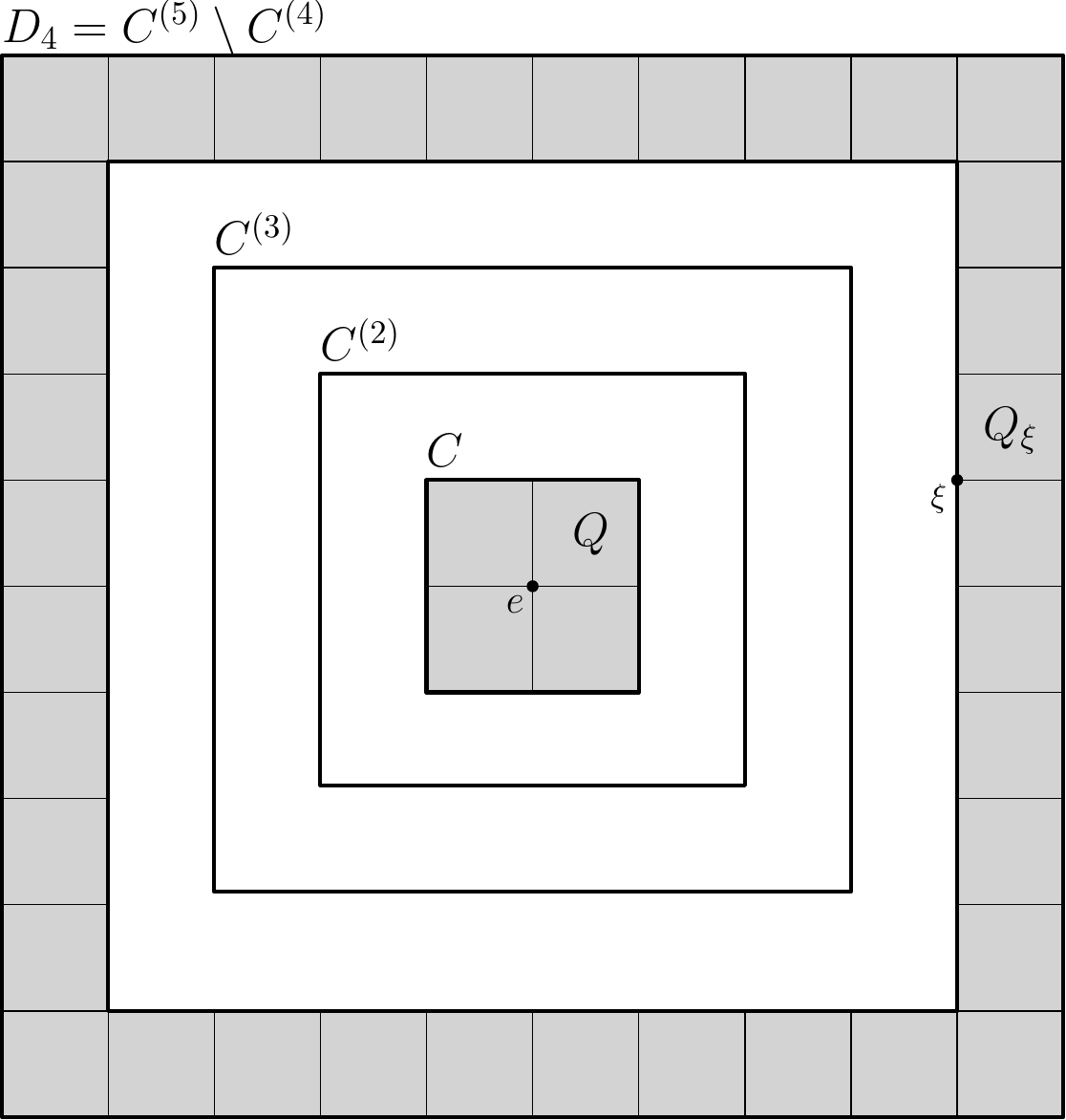}
           \par\vspace{0cm}
           \caption{Scheme of $D_4=C^{(5)}\setminus C^{(4)}=\bigcup_{\xi\in J_{5}\setminus J_4} Q_\xi $}
           \label{dactilios}
\end{figure}

\bdem
Since $\G$ is compactly generated, it can be represented by \eqref{cg} and  
$$
\G = \bigcup_{\xi \in \Xi} Q_{\xi},
$$
where $\Xi$ is the lattice given by \eqref{Ksi} and $Q_{\xi}$ is given by \eqref{Qtau}.
Define $$A_{\xi}:= \sup_{x\in Q_{\xi}} |\hat{\psi}(x)|,$$ then we have that
\begin{align*}
 \sum_{\gamma\in \Gamma: \gamma \neq e}|\hat{\fii}\big(r\,q_1(\gamma)\big)\hat{\psi}\big(q_2(\gamma)\big)|&=
\sum_{\xi \in \Xi} \sum_{\,q_2(\gamma) \in\, Q_{\xi}: \gamma \neq e}|\hat{\fii}\big(r\,q_1(\gamma)\big)\hat{\psi}\big(q_2(\gamma)\big)|\\
&\leq\sum_{\xi \in \Xi}A_{\xi} \sum_{\,q_2(\gamma) \in Q_{\xi}: \gamma \neq e}|\hat{\fii}\big(r\,q_1(\gamma)\big)|.
 \end{align*}
First, we will show that  for every $\eps>0$, there exists a finite set $J \subset \Xi$ such that
\begin{equation}
\sum_{\xi  \notin  J} A_{\xi} <\eps.
\end{equation}
Indeed, since $\hat{\psi}  \in  \mathcal{S} (\G)$  by the Definition \ref{A-functions}
we have that there exists $C=C_{\hat{\psi}}$ 
such that 
\begin{equation} \label{acot}
\| \hat{\psi}\chi_{\G  \setminus  C^{(k)}} \|_\infty\leq B_r k^{-r},
\end{equation}
where  $C^{(2)}=C+C$ and $C^{(k)}=C^{(k-1)}+C$ and $B_r$ a constant depending on $r$. Without loss of generality, we consider 
$$
C= [-n,n]^d \times \{-m,\ldots, m\}^{\ell} \times \K.
$$
Let  $D_k = C^{(k+1)} \setminus C^{(k)}$
and define $A_{D_k}:= \sup_{x\in D_k} |\hat{\psi}(x)|$. Fix $N>0$ to be chosen later on,  and let $J_N \subset \Xi$ such that 
$C^{(N)}= \bigcup_{\xi \in J_N} Q_{\xi}$ (see Picture \ref{dactilios}). Then, 
$$
\sum_{\xi \notin J_N} A_{\xi} \leq \sum_{k=N+1}^{\infty} A_{D_k} (2kn)^d (2km)^{\ell} \leq (2n)^d(2m)^{\ell} B_r \sum_{N+1}^{\infty} k^{-r+d+l},
$$
where in the last inequality we have used \eqref{acot}.  Choosing $r> d+\ell+1$ and $N$ big enough we have that $\sum_{\xi \notin J_N} A_{\xi} <\eps$, 
for every $\eps>0$. So, we can take $J$ equal to such $J_N$. 

\medskip

By Lemma \ref{separacion uniforme}, the quasicrystals $M_{Q_{\xi}}$ are uniformly bounded below by certain $\delta>0$.
Then, since  $\hat{\fii}\in\ese(\R)$ we have
$$
\sum_{\gamma \neq e:\,q_2(\gamma)\,\in\, Q_{\xi}}|\hat{\fii}\big(r\,q_1(\gamma )\big)|\leq 2^m
 \|\hat{\fii}\|_\infty+ \sum_{x \in M_{Q_{\xi}},\|x\|_\infty \geq \delta} \frac{1}{(1+\|r x\|_{\R^m}^2)^m}\ .
$$
Consequently, for $r$ big enough we have that 
\begin{align*}
\sum_{\xi \in \Xi \setminus J_N} 
A_{\xi} \sum_{\gamma:\,q_2(\gamma)\,\in\, Q_{\xi}}|\hat{\fii}\big(r\,q_1(\gamma)\big)|<(2^m \|\hat{\fii}\|_\infty+1)\eps. 
 \end{align*}
On the other hand, if $K=\bigcup_{\xi \in J_N} Q_{\xi}$,  using that $M_K$  is uniformly discrete (Lemma \ref{sep})  and so the points of $M_K$  are uniformly discrete by the identity element, we have
\begin{align*}
\sum_{\xi \in J_N}A_{\xi} \sum_{\gamma \neq e:\,q_2(\gamma)\,\in\, Q_{\xi}}|
\hat{\fii}\big(r\,q_1(\gamma)\big)|&\leq \|\hat{\psi}\|_\infty \sum_{x \in M_K,\, x  \neq 0} |\hat{\fii}(rx)|\\
&\leq C'\|\hat{\psi}\|_\infty \sum_{x \in M_K,\, x \neq 0} \frac{1}{(1+\|r x\|_{\R^m}^2)^m}\\
&=o(r^{-1}). 
 \end{align*}

  \edem

\bdem [Proof of Theorem \ref{nl}]
Using  Poisson's summation formula for $x=e$ we get 
\begin{align*}
\sum_{h\in H}\frac{1}{r^m}\fii\left(\frac{p_1(h)-a}{r}\right)\psi(p_2(h))&=\frac{1}{|H|} \sum_{\gamma \in\Gamma }e^{-2\pi i a\cdot q_1(\gamma)}\hat{\fii}\big(r\,q_1(\gamma)\big)\hat{\psi}\big(q_2(\gamma)\big).
\end{align*}
Then, as $\hat{\fii}\in\ese(\R^m) $ and $\hat{\psi}\in\ese(\G)$, by Lemma \ref{limit}
\begin{align}
 \sum_{\gamma \in\Gamma,\,q_2(\gamma)\neq e}|\hat{\fii}\big(r\,q_1(\gamma)\big)\hat{\psi}\big(q_2(\gamma)\big)|\xrightarrow[r\to\infty]{}0\label{que pollo}.
\end{align}
Since $p_1(H)$ is dense in $\R^m$, the projection $\left. q_2\right|_\Gamma$ is one to one. Therefore, $q_2(\gamma)=e$ if and only if $\gamma=(0,e)\in \R^m\times \G$. In consequence, only the term corresponding to $\gamma=(0,e)$ does not appear in \eqref{que pollo}, which is precisely $\hat{\fii}(0)\,\hat{\psi}(e)$.
\edem

\subsection{Beurling type theorem for simple quasicrystals} \label{BeurlingQ}

Let $(\R,G,H)$ be an complete CP-scheme and  $I \subseteq \R$ an interval.  Consider the associated  simple quasicrystal:

\begin{equation}
\L_I :=\{p_2(h): h \in H, \,p_1(h)\in I \} \subseteq G.
\end{equation}

\medskip

Before  proceeding to the proof of Theorem \ref{thsuf}, we need  an estimate of the Beurling density of $\L_I$.
The result is due to Meyer (see \cite{Meyer5}). If $G=\R^n$, the reader is also referred to  \cite{Meyer2}.
We provide here an alternative proof, based in the duality theorem and the extension of Landau's theorem 
to LCA groups.  Recall that for a compact set $K\subseteq \G$,
\begin{equation*} 
M_K := \{q_1(\gamma): \gamma\in \Gamma, \, q_2(\gamma) \in K \} \subseteq \R,
\end{equation*}
where  $\Gamma$ denotes the dual lattice  of  $H$.

\begin{teo} \label{un}
Given an interval $I\subseteq \R$,  
$$
\cD(\La_I) = \frac{|I|}{\m{H}}.
$$
\end{teo}
\bdem
Given $\eps>0$, consider a Riemann integrable compact set $C\subset \R$, such that 
$$
\mu_{\G}(C)= \frac{|I|+\eps}{\m{H}}.
$$
Taking into account Corollary \ref{uno} we have that $\cD(M_C)=|I|+\eps$. Then, since $M_C\subset \R$, we use Beurling Theorem \ref{TeoBeurling}, and so we get that $M_C$ is a stable sampling set for $PW_{\widetilde{I}}$, where ${\widetilde{I}}$ is a slight dilation of $I$ so that $|I|+\eps >{|\widetilde{I}|}$. 
 By the duality theorem, this implies that $\La_{\widetilde{I}}$ is an stable interpolation set for $PW_C$. So, by the extension to LCA groups of Theorem 
 \ref{LCA Landau}, we get that 
$$
\cD^+(\La_I)\leq\cD^+(\La_{\widetilde{I}})\leq \mu_{\G}(C)=\frac{(|I|+\eps)}{\m{H}}.
$$
 Since $\eps>0$ is arbitrary, we get that $\displaystyle \cD^+(\La_I)\leq  \frac{|I|}{\m{H}}$. In a similar way, using the other part of the duality theorem, 
 we can prove that $\displaystyle \cD^-(\La_I)\geq \frac{|I|}{\m{H}}$, which completes the proof.
\edem

Combining the result in Theorem \ref{un} with the duality theorem, we directly obtain Theorem \ref{thsuf} an extension of Beurling's theorem for simple quasicrystals,

\bdem[Proof of Theorem \ref{thsuf}]
Assume that  $\cD(\L_I) >  \mu_{\G}(K).$ 
Then by Theorem \ref{un},  
$$
 |I| > \frac{\mu_{\G}(K)}{|H|^{-1}}=  \frac{\mu_{\G}(K)}{\Gamma}= \cD(M_K).
 $$
Now we use Beurling Lemma in $\R$ (Theorem \ref{TeoBeurling}) to conclude that $M_K \subseteq \R$ is a set of stable interpolation for $PW_I.$
Therefore by the duality theorem $\L_I$ is a set of stable sampling for $PW_K$. 
The second claim of the theorem can be proved analogously.
\end{proof}

\begin{rem}
In \cite{neclca} Gr\"ochenig,  Kutyniok, and  Seip raised the question on the existence of stable sampling sets and stable interpolation set for $PW_{\Omega}$, whose densities were arbitrarily close to the critical density given by (the generalized) Landau's theorem. Although this question was completely solved in \cite{AAC} using another approach,  it is worth to mention that Theorem \ref{thsuf} gives an alternative proof of the existence of these sets, providing the group $G$ admits a quasicrystal. 
\end{rem}

%-----------------------------------------------------------------------------------------------------------------------
%-----------------------------------------------------------------------------------------------------------------------
% APPENDIX
%-----------------------------------------------------------------------------------------------------------------------
%-----------------------------------------------------------------------------------------------------------------------

\appendix{
\section{Divisible groups}\label{App}
A (discrete) abelian group $\D$ is called \textbf{divisible} is for every $\si\in\D$ and $n\in \N$ there exists $b\in\D$ such that
$$
nb=\si.
$$
The easiest example of a divisible group is the group of rational numbers with the addition $(\Q,+)$. Another example of divisible groups are the so called quasi-cyclic groups $\Z(p^\infty)$, also called Pr\"ufer groups. Recall that, given a prime number $p$, then 
$$
\Z( p^{\infty} ) := \left\{r\in [0,1): r=\frac{k}{p^n}; \, k\in \N_0\,,\ n\in \N \right\}, 
$$
endowed with the usual addition mod one. Although these are not all the divisible groups, they are the building blocks of any countable 
divisible group  (see \cite{F}, pg. 104).
\begin{teo} \label{divisible}
Any countable divisible group $\D$ is of the form  
$$
\D\simeq \Q^{k_q}  \oplus   \bigoplus_{p: \text{ prime} }\Z(p^{\infty})^{k_p},
$$
where $k_q$ and all the $k_p$ are non-negative integers, which constitute a complete invariant of $\D$.
\end{teo}
 
\medskip
 
Our interest in divisible groups comes from the following universal property  (see \cite{F}, pg. 106):

\begin{teo}
Any countable group can be embedded as a subgroup of a divisible group. 
\end{teo} 

\medskip

Let us identify  $\Z_p$ with the group $\{k/p: 0\leq k <p\}$ endowed with the addition mod 1. Note that if a discrete group $\D$ has a copy of $\Z_r^n$ for some prime number $r$, and $\D$ is embedded in a divisible group
\begin{equation}\label{eq divisible}
\Q^{k_q}  \oplus   \bigoplus_{p: \text{ prime} }\Z(p^{\infty})^{k_p},
\end{equation}
then the value of $k_r$ has to be greater or equal than $n$. Indeed, the only components in \eqref{divisible} 
that have elements of order $r$ are the $k_r$ copies of $\Z(r^{\infty}).$  This observation leads to the following result.

\begin{pro}\label{como son}
Let $\D$ be a  countable abelian group such that, for any prime number $p$, it does not have a copy of $\Z_p^{m+d+1}$. Then it is isomorphic to a subgroup of a divisible group of the form
$$
\Q^{k_q}  \oplus   \bigoplus_{p: \text{ prime} }\Z(p^{\infty})^{k_p},
$$
such that $k_p\leq m+d$ for every prime $p$. 
\end{pro}

}
\vspace{.5cm}

\noindent
\begin{center}
\small{ACKNOWLEDGEMENT}
\end{center}
We would like to thank Yves Meyer for suggesting us to study the relevant problem of the existence of quasicrystals in LCA groups, that generated the results in Section~\ref{existence of lattice}. We would also like to thank the referee for his/her careful reading and suggestions, which help us to make the article more readable.

%------------------------------------------------------------------------------
%------------------------------------------------------------------------------
\end{document}